\documentclass[11pt]{amsart}

   \usepackage[a4paper, margin=100pt]{geometry}    
    \usepackage[T1]{fontenc}
    \usepackage[british]{babel}
    \usepackage[utf8]{inputenc}

    \usepackage{mathabx}

\normalsize 
    \addtocontents{toc}{\setcounter{tocdepth}{1}}

    \usepackage{xcolor}
    \usepackage{hyperref}
         \hypersetup{colorlinks, citecolor=blue, filecolor=blue, linkcolor=black, urlcolor=blue}
    \usepackage{soul}
    \usepackage[numbers]{natbib}    
   
    \usepackage{caption}
    \usepackage{subcaption}

    \usepackage{tikz-cd}

\usepackage{marginnote}

\newcommand{\exwar}[1]{}
\newcommand{\ua}{\underline{a}}
\newcommand{\ub}{\underline{b}}
\newcommand{\uc}{\underline{c}}
\newcommand{\uw}{\underline{w}}

\newcommand{\ik}{\ldbrack -k, k\rdbrack}

    \newcommand{\sldr}{\mathrm{SL}(d,\mathbb{R})}

    \newcommand{\transition}[1]{[{#1}\rangle}

\DeclareMathOperator{\GL}{\rm{GL}}
\DeclareMathOperator{\SL}{\rm{SL}}

\DeclareMathOperator{\cl}{cl}

\newcommand{\R}{\mathbb{R}}\newcommand{\N}{\mathbb{N}}
\newcommand{\Z}{\mathbb{Z}}

\newcommand{\id}{\mathrm{Id}}

\newcommand{\cA}{\mathcal{A}}
\newcommand{\cB}{\mathcal{B}}
\newcommand{\cC}{\mathcal{C}}
\newcommand{\cD}{\mathcal{D}}

\newcommand{\cF}{\mathcal{F}}
\newcommand{\cG}{\mathcal{G}}
\newcommand{\cH}{\mathcal{H}}

\newcommand{\cO}{\mathcal{O}}

\newcommand{\cQ}{\mathcal{Q}}

\newcommand{\cS}{\mathcal{S}}

\newcommand{\cU}{\mathcal{U}}
\newcommand{\cV}{\mathcal{V}}

\newcommand{\bN}{\mathbb{N}}

\numberwithin{equation}{section}
\newtheorem{theorem}{Theorem}[section] %[section]
\newtheorem{corollary}[theorem]{Corollary}
\newtheorem{lemma}[theorem]{Lemma}
\newtheorem*{lemma*}{Lemma}
\newtheorem{proposition}[theorem]{Proposition}
\newtheorem*{proposition*}{Proposition}

\newtheorem{question}[theorem]{Question}

\newtheorem*{question*}{Question}
\newtheorem*{theorem*}{Theorem}
\newtheorem*{claim*}{Claim}

\newtheorem{theoremain}{Theorem}

\theoremstyle{definition}
\newtheorem{definition}[theorem]{Definition}

\newtheorem{example}[theorem]{Example}

\theoremstyle{remark}
\newtheorem{remark}[theorem]{Remark}

    \usepackage{bigints}
    \usepackage{xfrac}
    \usepackage{mathrsfs}
    \usepackage{comment}
    \usepackage{verbatim}
    \usepackage{adjustbox}
    \usepackage{graphicx}
    \usepackage{import}

   \newcommand{\gldr}{\mathrm{GL}(d,\R)}
    \newcommand{\cylinder}{\mathsf{C}}
    \newcommand{\qc}{\cQ\cC}

\begin{document}

    \title[cocycles with quasi-conformality]{Cocycles with Quasi-Conformality I:\\ Stability and abundance}

\author[M. Nassiri]{Meysam Nassiri}
\address{School of Mathematics, Institute for Research in Fundamental Sciences (IPM), P.O. Box 19395-5746, Tehran, Iran}
\email{nassiri@ipm.ir}
\author[H. Rajabzadeh]{Hesam Rajabzadeh}
\address{School of Mathematics, Institute for Research in Fundamental Sciences (IPM), P.O. Box 19395-5746, Tehran, Iran}
\email{rajabzadeh@ipm.ir}
\author[Z. Reshadat]{Zahra Reshadat}
\address{School of Mathematics, Institute for Research in Fundamental Sciences (IPM), P.O. Box 19395-5746, Tehran, Iran}
\email{zahrareshadat@ipm.ir}

\thanks{
 This work is partially supported by INSF grant no.  4001845.\\ \indent Z.R. was supported by the IMU Breakout Graduate Fellowships.}
\date{}

\begin{abstract}
This is the first part of a series of papers devoted to the study of linear cocycles over chaotic systems. In the present paper, we establish the existence of such cocycles that $\cC^\alpha$-stably exhibit fiberwise bounded orbits ($\alpha>0$). The proof is based on a new mechanism which yields stable elliptic-type behavior in $\gldr$ or $\sldr$ cocycles. Moreover, we show that this phenomenon is $\cC^0$-dense among
$\sldr$ cocycles over a shift of finite type without dominated splitting.

\end{abstract}

\maketitle
\tableofcontents

\section{Introduction}\label{sec:intro}

A main goal in the field of dynamical systems 
is unveiling stable properties that cannot be destroyed by small perturbations. Another fundamental goal is to identify typical properties among all systems or within specific classes.
In this paper, we address both of these objectives for linear cocycles.

Given a homeomorphism $f:X\to X$ of a compact metric space $X$, and a vector bundle $\pi:\boldsymbol{E} \to X$,
a \textit{linear cocycle} over $f$ is a homeomorphism $F:\boldsymbol{E} \to \boldsymbol{E}$ such that $\pi \circ F=f\circ \pi$, and the restriction of $F$ to the fibers of $\boldsymbol{E}$ is linear. 
A canonical example is the derivative $Df:TM \to TM$ of a diffeomorphism $f$ on a manifold $M$. 
Another illuminating example that naturally leads to a cocycle defined on a trivial vector bundle is the cocycle associated with the random products of a finite family of matrices.

The principal subject of study in linear cocycles is the behavior of fiberwise dynamics, where concepts such as hyperbolicity, dominated splitting, and Lyapunov exponents are developed. This study is at the core of many branches in dynamical systems and related fields, including smooth ergodic theory, homogeneous dynamics, and mathematical physics. 

Uniform hyperbolicity and its extensions, such as dominated splitting, are stable phenomena that give rise to rich and chaotic dynamics (see \S \ref{subsec:def-domination} for the definitions in the case of linear cocycles). 
As a matter of fact, there are open sets of systems that do not admit dominated splitting. 
Therefore, it is fundamental to identify stable and prevalent behaviors appearing in the lack of domination.
\begin{question}\label{ques:stable-no-domination}
    What stable phenomena prevent a linear cocycle from admitting (global) dominated splitting? 
\end{question} 
A consequence of dominated splitting for a cocycle is a uniform gap in the Lyapunov spectrum with respect to every measure invariant under the base dynamics. Hence, an obstruction for global dominated splitting is the degeneracy of the Lyapunov spectrum for some invariant measure on the base. Several works have employed this approach to describe the cocycles without domination. 
In this direction, we can mention the statements announced by Mañé in ICM 1983 \cite{Mane-ICM,Mane-Announced}, later proved by Bochi {\cite{bochi2002,Bochi-2010}} and Bochi-Viana \cite{Bochi-Viana2002,Bochi-Viana}, and improved  by Avila-Crovisier-Wilkinson \cite{ACW-ihes}. 
All these results provide a dichotomy between a degenerate spectrum and dominated splitting for generic linear cocycles. 
More precisely, they show that for systems in a residual subset (in certain classes) of volume-preserving $\cC^1$ diffeomorphisms of closed smooth manifolds, either the derivative cocycle admits a variant of hyperbolicity or all the  Lyapunov exponents for almost every point vanish. The statements also apply in the setting of continuous linear cocycles over certain bases. There are also other works aimed to characterize linear cocycles without domination in terms of the Lyapunov spectrum under various conditions on the cocycles, namely on the base dynamics or on the dimension of the fibers. Among them we can mention  \cite{Bochi-Fayad,Bochi-Navas,Avila-Bochi,Avila-Bochi-TAMS,Bochi-2010,Bochi2013,Potrie-Kassel}. See also the survey paper  \cite{bochi-vianasurvey} and the book \cite{Viana_2014} and references therein. 

The aforementioned works describe the generic behavior of fiberwise dynamics of linear cocycles. Nevertheless, they do not provide stable phenomena in the absence of dominated splitting. 
{Indeed, concerning Question \ref{ques:stable-no-domination}, one needs a mechanism ensuring the persistence of ergodic invariant measures with degenerate fiberwise spectrum under small perturbations. This is the strategy followed in \cite{BBD-robust-vanishing-all} to identify open sets of iterated function systems with vanishing spectrum (see also \cite{BBD-2016,BBD-2018,Diaz-ICM,BDK2021}).}
{Another strategy to achieve stable behaviors obstructing dominated splitting is based on focusing on the periodic measures (orbits) and highlights the notion of ellipticity.}
If a periodic fiber exists such that the eigenvalues of the composition of associated linear maps up to the period have the same modulus, then the cocycle cannot admit global dominated splitting. For volume-preserving cocycles (those that preserve the volume form on the fibers), this condition is equivalent to the existence of a periodic fiber that is conjugated to an isometry.

In volume-preserving cocycles over 2-dimensional vector bundles, it follows from the work of Mañé \cite{Mane} that if the base dynamics is sufficiently chaotic, a generic dichotomy of elliptic/hyperbolic behavior holds; a cocycle in an open and dense subset either admits uniform hyperbolicity or admits a periodic point with non-real eigenvalues of modulus 1. 
An analogous result for linear iterated function systems generated by a finite family of elements of $\SL(2,\R)$ has been proved in \cite{Yoccoz2004,Avila_Bochi_Yoccoz}, where the authors give a delicate understanding of the set of generators, leading to ellipticity or hyperbolicity.

In attempting to extend such dichotomy to higher dimensions, the initial observation is that unlike dimension 2, the existence of periodic fibers with eigenvalues of equal modulus is not stable in higher dimensions. 
Nevertheless, the work of Bonatti-Diaz-Pujals \cite{BDP} shows that every (volume-preserving) linear cocycle over certain chaotic systems (such as the shift map) either admits dominated splitting or has several periodic points with a pair of complex conjugate non-real eigenvalues (cf. \cite{Gourmelon2022}).
They use such periodic points to show that every volume-preserving cocycle without dominated splitting
can be $\cC^0$-approximated by cocycles with identity action on some periodic fiber (cocycles with homothetic action on some periodic fiber if the cocycle is not volume-preserving). See Section \ref{sec:transition} for more details.
So, when the base dynamics $f$ and the fiber bundle $\boldsymbol{E}$ are fixed, if we denote the set of continuous cocycles $F$ admitting dominated splitting by $\cD$, and the cocycles which are conjugated to an isometry on some periodic fiber by $\cB_{\rm Per}$, the equality $\cD^c = \cl(\cB_{\rm Per})$ holds in the space of volume-preserving cocycles endowed with the $C^0$ topology.
As a matter of fact, $\cB_{\rm Per}$ has no interior, and its complement is residual for $d>2$.

In this paper, we study the existence of points (not necessarily periodic) whose associated linear maps under the iterations of cocycles remain (uniformly) close to homothety.

 More precisely, we say 
 that the linear cocycle $F$ over $f:X\to X$ \emph{admits fiberwise bounded orbits} if there exists $x\in X$ and $\kappa>0$ such that 
\[\|F_n(x)^{\pm 1}\|<\kappa,\quad \text{for every}~ n\in \Z.\]
where $F_n(x)$ denotes the linear map induced from $F^n$, the $n$-th iteration of $F$ to fiber $\boldsymbol{E}(x):=\pi^{-1}\{x\}$.
Similarly, we say that $F$
\emph{admits quasi-conformal orbits}, if there exists $x\in X$ and $\kappa>1$ such that  
\[\|F_n(x)\|.\|F_n(x)^{-1}\|<\kappa, \quad \text{for every}~ n\in \Z,\]
Whenever the fiber bundle $\boldsymbol{E}$ is a trivial vector bundle, the linear maps $F_n(x)$ correspond to some matrices. The above notions are equivalent to saying that the associated product of matrices remains in a specific set.

If $\boldsymbol{E}$ is a trivial vector bundle
of dimension $d\in \N$, every linear cocycle $F$ is represented by a continuous function $A:X\to \gldr$; that is for every $x\in X$ and $v\in \R^d$, $F(x,v)=(f(x),A(x)v)$. In such a case, we write $(f,A)$ instead of $F$ and say that $(f,A)$ is a $\gldr$ cocycle. The concept of $\sldr$ cocycles is defined similarly.  Clearly, for every $n\in \Z$, $F^n$ is a cocycle over $f^n$. We denote the matrix-valued function associated with $F^n$ by $A_n$. The cocycle $(f,A)$ admits fiberwise bounded (resp. quasi-conformal) orbits if there exists an open neighborhood $\cU$ of $\id$ (resp. of orthogonal matrices) in $\gldr$ with compact closure such that $A_n(x)\in \cU$, for every $n\in \Z$.

We denote the set of cocycles admitting quasi-conformal orbits by $\qc$, and the set of cocycles admitting fiberwise bounded orbits by $\cB$.  These two notions are closely related. They coincide in volume-preserving cocycles. For general linear cocycles, by factoring out the determinant, cocycles in $\cB$ will correspond to cocycles in $\qc$. Note that, $\cB\subset \qc$, but $\cB\neq \qc$. Indeed, one can start with a cocycle in $\qc$, for instance, by taking all fiberwise linear maps to be {isometry (which is contained evidently in both $\cB,\qc$) and multiply all the fiberwise linear maps with a constant number greater than 1. The resulting cocycle is again in $\qc$ but not in $\cB$.

On the other hand, the set  $\qc$ (and hence $\cB$) does not intersect $\cD$. The aim of this paper is to understand the properties of sets $\qc$ and $\cB$. 
A question in this direction is the following. 
\begin{question}\label{ques:interior-B2}
When the base dynamics $f$ is hyperbolic, 
    does $\cB$ have a non-empty interior?
\end{question}

Our first result answers this question affirmatively for cocycles over the shift map in the $\cC^\alpha$ topology.

\begin{theoremain}\label{thm:GL-bd-example}
For every $\alpha>0$ and every integer $d\geq 2$, in the space of $\gldr$ cocycles over {the full shift on finitely many symbols}, the $\cC^\alpha$-interior of $\cB$ is non-empty.
\end{theoremain}

The definition of Hölder topology on the space of linear cocycles is provided in \S\ref{sec:holder}.  
{We would like to remark that the Hölder topology is essential in our proof.
Indeed, very recently Bochi \cite{Bochi-no-stable-CQ} has shown that the $\cC^0$-interior of $\cB$ is empty. Further discussion can be found in Subsection \ref{subsec:question-C0-stable}.}

The next question concerns 
prevalence of cocycles in $\cD^c$ with fiberwise bounded orbits. Clearly, the mentioned example made by perturbation of a cocycle with orthogonal elements demonstrates that $\cB$ is not dense even in the complement of $\cD$. However, one may expect the density of cocycles (stably) admitting quasi-conformal orbits in $\cD^c$.

\begin{question}\label{ques:interior-B}
When the base dynamics $f$ is hyperbolic, 
is the {$\cC^\alpha$-}interior of $\qc$  dense in $\cD^c$? 
\end{question}

{The answer to this question may depend on the topology on the space of cocycles. Here, in this paper, we address the density in Question \ref{ques:interior-B}  with respect to the $\cC^0$ topology.} {Note that for volume-preserving cocycles since $\cB_{\rm Per}\subset \cB$, it follows from the work of Bonatti-Diaz-Pujals \cite{BDP} that the set  $\cB$ is dense in $\cD^c$. Nevertheless, this does not imply the density of the interior of $\cB$ in $\cD^c$ when the dimension of fibers is greater than 2 as $\cB_{\rm Per}$ has no interior. Indeed, in this case, for a generic cocycle (in the $\cC^\alpha$ topology), the linear maps on all periodic fibers on their period are hyperbolic (see Proposition \ref{prop:generic hyp periodic}).
}

{Here, we use the method introduced in \cite{FNR} based on a covering condition to establish orbits whose fiberwise linear maps remain in a bounded region along iterations of the cocycle, and we can guarantee this covering condition, sufficiently $\cC^0$-close to every cocycle. As a result, we give a positive answer to Question \ref{ques:interior-B} in Hölder regularity.}

\begin{theoremain}\label{thm:GL-QC}
  For every $\alpha>0$ and every integer $d\geq 2$, in the space of $\gldr$ cocycles over {the full shift on finitely many symbols}, 
 the $\cC^\alpha$-interior of $\qc$ is  $\cC^0$-dense in $\cD^c$.
\end{theoremain}

Consider the shift map $\sigma$ acting on the space of bi-infinite sequences from a finite alphabet, and let $(\sigma,A)$ be a $\gldr$ cocycle over $\sigma$. The continuous $\gldr$ cocycle  $(\sigma,A)$ admits quasi-conformal orbits if and only if the cocycle $(\sigma,\hat{A})$ defined by $\hat{A}(x)=|\det A(x)|^{-1/d}A(x)$ admits fiberwise bounded orbit. This is because $|\det \hat{A}(x)|=1$. 
So Theorem \ref{thm:GL-QC}  follows from a similar theorem for volume-preserving cocycles. 
The following is the analog of Theorem  \ref{thm:GL-QC} for $\sldr$ cocycles. 

\begin{theoremain}\label{thm:SL-B}  For every $\alpha>0$ and every integer $d\geq 2$, in the space of $\sldr$ cocycles over {the full shift on finitely many symbols}, 
    the $\cC^\alpha$-interior of $\cB$ is  $\cC^0$-dense in $\cD^c$. 
\end{theoremain}

Our mechanism provides a $\cC^\alpha$ open set of cocycles with bounded fiberwise orbits on non-periodic fibers. This is a new result even for 2-dimensional cocycles. 
{
Our approach relies on the fact that locally constant cocycles are $\mathcal{C}^0$-dense in the space of all continuous cocycles over the shift. Indeed, the proofs of our main theorems imply the following statement for locally constant cocycles.}

\begin{theoremain}\label{thm:lc-dichotomy}
Every locally constant $\mathrm{SL}(d,\mathbb{R})$ cocycle over {the full shift on finitely many symbols} either admits a dominated splitting or can be $\mathcal{C}^0$-approximated by locally constant cocycles that $\mathcal{C}^\alpha$-stably admit fiberwise bounded orbits.
\end{theoremain}

{
Although Theorem \ref{thm:lc-dichotomy} has some apparent similarity to some results in \cite{Avila_Bochi_Yoccoz,Yoccoz2004} on the generic behavior of iterated function systems generated by a finite family of matrices in $\mathrm{SL}(2,\R)$, there are substantial differences. 
 In those works, all the cocycles are restricted to the space of cocycles that depend only on the zeroth position. In contrast, in our setting, the locally constant cocycles do not necessarily depend only on the zeroth position, and the perturbations involve cocycles depending on a larger number of positions. Therefore, the proofs in the present paper do not address analogues of Theorem \ref{thm:SL-B} within the space of locally constant cocycles depending only on the zeroth position; see Subsection \ref{ques:generic-IFS} for further discussion.
}

We remark that, {thanks to the symplectic variants of \cite{BDP} established in \cite{Horita-Tahzibi}}, the proofs of Theorems \ref{thm:SL-B} and \ref{thm:lc-dichotomy} can be modified to apply also to cocycles with values in $\mathrm{Sp}(2d, \mathbb{R})$, the group of symplectic $2d \times 2d$ matrices. See \cite{Reshadat-thesis} for a detailed proof.

\subsection{How big is the set of fiberwise bounded orbits?}

{
In order to analyze the quasi-conformal orbits for a general linear cocycle $(\sigma,A)$ over the shift map, we may assume that the map $A$ takes values in $\mathrm{SL}(d, \mathbb{R})$. A key question regarding quasi-conformal behavior on the fibers is the prevalence of these orbits. To better understand the set of fiberwise bounded orbits, it is natural to identify the invariant measures it supports. From the viewpoint of Bernoulli measures on the shift space, this set has measure zero for most cocycles. Indeed, for a majority of $\alpha$-Hölder maps $A$, the cocycle $(\sigma,A)$  
has a positive top Lyapunov exponent (see \cite{Furstenberg,Ledrappier,VianaNonvanishing_of_exponents,invariance-principle}). In the forthcoming paper \cite{NRR2}, we will provide a positive answer to the problem of \emph{the existence of ergodic measures with positive entropy supported on the set of bounded orbits}, shedding light on the topology of the set of bounded orbits. This is relevant to a series of works concerning the robust existence of non-hyperbolic ergodic invariant measures with positive entropy (see \cite{BDK2021,Gelfert-Rams,Lacka}).
}

\subsection{Critical behaviors of dynamical systems}

Critical behaviors of dynamical systems usually refer to phenomena originated to (stable) non-hyperbolicity, mainly because the  nearby systems may exhibit different dynamical behaviors. 

For one-dimensional smooth maps, critical points may also refer to points with vanishing derivative. A remarkable theorem of Ma{\~n}{\'e} \cite{Mane85} shows that a generic one-dimensional smooth map without
any point with vanishing derivative is either hyperbolic or conjugate to irrational rotation.

For diffeomorphisms in dimension two, homoclinic tangencies
play a role similar to the critical points of one dimensional maps (see  \cite{BC,MV} in the study of  H{\'e}non maps). {The works of \cite{PR} and \cite{CP} developed the concept of critical points for surface diffeomorphisms, pursuing the idea of defining the critical set \emph{as the region where domination fails}}. This resembles the results of Ma{\~n}{\'e} in dimension one, and as expected, the orbit of a homoclinic tangency is an essential example of a critical orbit. 
However, the following fundamental question (see \cite[\S 4.2.2]{BDV}) has remained a big challenge particularly in higher dimensions, despite the rich developments concerning homoclinic bifurcations over the last five decades (See \cite{PT,Palis-conj-2005,BDV,BD-2012,BBD-2016} and references therein). 
\begin{question}
   Is there an intrinsic definition of critical points or critical orbits for diffeomorphisms? The same question may be asked for linear cocycles, group actions, etc. 
\end{question}

The results of this paper suggest a candidate for {\it super-critical} orbits in linear cocycles over shift map (or other chaotic base dynamics), namely, ``points with fiberwise bounded orbit''. Also, we can consider points with uniformly bounded derivatives in the family of smooth cocycles.

With this idea in mind, we expect that points with fiberwise bounded (or quasi-conformal) orbits in the smooth action of groups will have consequences similar to the role played by critical orbits in the statistical behavior of typical dynamical systems (See \cite{Doug-Melbourne-Talebi,Doug-Luzzatto}).
An answer to the next question may contribute to a better understanding of the origins of the so-called historical behavior and may shed some light on Palis' program \cite{Palis-conj-1995, Palis-conj-2005}.

\begin{question}
    What are the statistical consequences of the existence of super-critical orbits? Does it yield non-statistical (or historic) behavior?
\end{question}

In symplectic dynamics, elliptic points are responsible for critical behavior. Elliptic points refer to periodic points with derivatives conjugated to an isometry with non-real eigenvalues.
On the other hand, for conservative dynamics in dimensions bigger than 2, ellipticity on a periodic fiber is unstable and rare. 
Indeed, typically all periodic fibers are hyperbolic  (See Lemmas \ref{lem:generic-pair-hyperbolic} and \ref{prop:generic hyp periodic}).
This raises an important question:

\begin{question}
   What are stable behaviors generalizing ellipticity in higher dimensions and beyond the symplectic setting?
\end{question}

It is convenient and fruitful to address this question for the derivative cocycle as well as for linear cocycles. 
Theorem \ref{thm:SL-B} may be seen as a response to this question for linear cocycles over the shift map. 
Analogous results for cocycles associated with smooth actions are also of interest. On the other hand, by considering the action of partially hyperbolic systems in the central bundle, the same notions can be defined. 
Such results and their implications for partially hyperbolic dynamics will be treated in a forthcoming paper.

\subsection{Ideas of the proofs and organization}
In Section \ref{sec:pre}, we present essential definitions and fix some notations that we will use in the subsequent sections. Section \ref{sec:qc and bounded} focuses on the basic properties of cocycles admitting quasi-conformal orbits or fiberwise bounded orbits. In the following sections, we aim to introduce a criterion for the existence of quasi-conformal orbits. It is based on a \emph{covering condition} for the action of a linear group on itself, initially introduced in \cite{FNR}.  There, it
has been used as the main tool in providing sequences of matrices with uniformly bounded products. Adapting this covering condition to the general locally constant cocycles requires the concept of \emph{transition}, inspired by the notion of transition in \cite{BDP}. A detailed discussion of this notion is given in \S \ref{subsec:loc-const-transition}.
Section \ref{Sec:covering} is devoted to developing and adapting this covering condition to the general locally constant cocycles over the shift map.  
  In Section \ref{sec:stability-qc}, we provide the proof of Theorem \ref{thm:GL-bd-example}. Indeed, we show that the covering condition for cocycles is sufficient to obtain a $\cC^\alpha$ open set of cocycles admitting quasi-conformal orbits.

  Next, in order to prove other theorems regarding the dichotomy of dominated splitting vs. stable quasi-conformality, we need to show that the covering condition can be satisfied for cocycles arbitrarily close to every cocycle without domination. 
   This is done in several steps. In Section \ref{sec:transition}, we recall and extend some results in \cite{BDP} to make periodic fibers with identity for the perturbation of volume-preserving cocycles in $\cD^c$. The extension lemma in Appendix \ref{sec:extension} allows us to extend a cocycle continuously after some perturbations at finitely many points. Finally, in Section \ref{sec:dichotomy}, we give the proofs of Theorem  \ref{thm:GL-QC}, \ref{thm:SL-B} and Theorem \ref{thm:lc-dichotomy}. The idea is to perturb a cocycle with many identity periodic fibers and get a locally constant cocycle satisfying the covering condition (for suitable open set and transitions). This, combined with the results of Section \ref{sec:stability-qc}, implies the existence of $\cC^\alpha$ open sets of cocycles admitting quasi-conformal orbits arbitrarily $\cC^0$ close to the initial cocycle.

\subsection*{Acknowledgement}
The authors would like to thank Jairo Bochi, Sylvain Crovisier, Abbas Fakhari, Enrique Pujals, and Amin Talebi for useful conversations and comments. Z.R. would also like to express her gratitude to Pierre Berger for providing an opportunity to visit the Jussieu Institute, where part of this paper was written. The authors are also grateful to the anonymous referee for valuable comments which  improved the exposition of the paper.

\section{Preliminaries} \label{sec:pre}
\subsection{Linear cocycle}  
Fix $d\in \N$. A $d$-dimensional \textit{real vector bundle} $(\boldsymbol{E},X,\pi)$ consists of topological spaces $X$ (base space) and $\boldsymbol{E}$ (total space), and a continuous  surjective map $\pi:\boldsymbol{E} \to X$ (bundle projection), where the following conditions are satisfied:
\begin{itemize}
    \item For every $x\in X$, $\boldsymbol{E}(x):=\pi^{-1}(\{x\})$ is a $d$-dimensional real vector space.
    \item (Local trivialization) For every $x\in X$, there exist an open neighborhood $U_x \subset X$ and a homeomorphism $\phi:U_x \times \R^d \mapsto \pi^{-1}(U_x)$, such that $\pi \circ\phi (x,v)=x$ and the map $v\mapsto \phi(x,v)$ is a linear isomorphism between the vector spaces $\R^d$ and $\pi^{-1}(\{x\})$.
\end{itemize}

The existence of local trivialization is equivalent to the existence of some orthonormal coordinate system for the fibers. 
{This is equivalent to being able to choose the local trivializations to be orthonormal. We assume that the fibers $\boldsymbol{E}(x)$ are endowed with inner products that vary continuously with the base point $x$. This induces a norm on each fiber, denoted by $\| \cdot \|_x$ (or simply $\| \cdot \|$), and a corresponding volume form. For notational simplicity, we will often refer to the vector bundle $(\boldsymbol{E}, X, \pi)$ simply as $\boldsymbol{E}$. 
}

Throughout the paper, we assume that the base space $X$ is a compact metric space and $f:X\to X$ is a homeomorphism. Also, we denote by $\cO(x)$ (resp. $\cO^+(x)$) the orbit (resp. the forward orbit) of $x$ under $f$. 

 By a \textit{linear cocycle}  over  $f:X\to X$, we mean a homeomorphism $F:\boldsymbol{E} \to \boldsymbol{E}$ such that $\pi \circ F=f\circ \pi$, and the restrictions of $F$ to the fibers of $\boldsymbol{E}$ are linear. We denote this cocycle by $(\boldsymbol{E},X,F,f)$ or simply by $F$. Also, the restriction of $F$ to the fiber $\boldsymbol{E}(x)$ is denoted by $F(x):\boldsymbol{E}(x)\to \boldsymbol{E}(f(x))$. 
{By choosing a linear  basis for fibers $\boldsymbol{E}(x)$, one can represent $F(x)$ by an element of $\gldr$. }

For  $n\in \Z$, the action of $n$-th iteration of this cocycle on the fiber $\boldsymbol{E}(x)$ is denoted by $F_n(x)$; thus for $n\in \N$
\begin{align*}
    F_n(x)&=F(f^{n-1}(x))\circ \cdots \circ F(x),\\
    F_{-n}(x)&=F(f^{-n}(x))^{-1}\circ \cdots \circ  F(f^{-1}(x))^{-1}.
\end{align*}
 Also, if $n=0$, $F_n(x)$ is the identity map on the fiber $\boldsymbol{E}(x)$.  The linear cocycle $F$ is \emph{volume preserving} whenever the maps $F(x)$
preserve the volume forms on the fibers. Equivalently,  by choosing an orthonormal family of bases for the fibers, all the matrices associated with linear maps $F(x)$ have determinant $\pm 1$.

 The space of all linear cocycles over $f$ endowed with the $\cC^0$ topology is denoted by  
 $\cC^0(\boldsymbol{E},f)$.  The $\cC^0$ topology on this space is induced by the distance 
\[ \textit{d}_{\mathcal{C}^0}(F,F'):= \sup_{x\in X} \| F(x)-F'(x) \|, \]
where $F,F'\in \cC^0(\boldsymbol{E},f)$ and $\|. \|$  denotes the operator norm of the fiberwise linear transformations. {Notice that $(\cC^0(\boldsymbol{E},f),\textit{d}_{\mathcal{C}^0})$ is a Baire space, that is, every countable intersection of dense open sets is dense.}
{For a subset $\cS \subset \cC^0(\boldsymbol{E},f)$, the closure of $\cS$ in $\cC^0$ topology is denoted by $\mathrm{cl}(\cS)$. When we do not specify the topology on the space of cocycles, we consider $\cC^0$ topology.}

\subsection{Dominated splitting}\label{subsec:def-domination}  We say that the vector bundle $\boldsymbol{E}'\to X$ is a vector subbundle of the vector bundle $\boldsymbol{E}\to X$, if for every $x\in X$, $\boldsymbol{E}'(x)$ is linear subspace of $\boldsymbol{E}(x)$.  In a similar fashion, we say the vector bundle $\boldsymbol{E}$ is the direct sum of  subbundles $\boldsymbol{E}_1,\boldsymbol{E}_2,\ldots,\boldsymbol{E}_k$ if for every $x\in X$, $\boldsymbol{E}(x)=\boldsymbol{E}_1(x)\oplus \cdots \oplus \boldsymbol{E}_k(x)$. 

For a cocycle $F:\boldsymbol{E}\to \boldsymbol{E}$ over $f$,  the above splitting is called $F$-invariant if $F(x)\boldsymbol{E}_i(x)=\boldsymbol{E}_i(f(x))$ for all $x\in X$ and $i=1,\ldots,k$.

\begin{definition}
    A splitting $\boldsymbol{E}=\boldsymbol{E}_1 \oplus \boldsymbol{E}_2$ into $F$-invariant non-trivial subbundles $\boldsymbol{E}_1,\boldsymbol{E}_2$ is called a {\it dominated splitting},  if there exists  $m\in \bN$ such that for every $x\in X$,
\[ \| F_m(x) \big|_{\boldsymbol{E}_2(x)} \| .\|{(F_m(x) \big|_{\boldsymbol{E}_1(x)})^{-1}} \| \leq \frac{1}{2}.\]
We use the notation $\boldsymbol{E}_2\prec \boldsymbol{E}_1$ if the above inequality holds. {For $j\in \N$, {\it dominated splitting of index $j$} means a dominated splitting where $\boldsymbol{E}_2$ is a $j$ dimensional bundle.}
\end{definition}

{In particular a cocycle $(\boldsymbol{E},X,F,f)$ is uniformly hyperbolic if there exist an $F$-invariant splitting $\boldsymbol{E}=\boldsymbol{E}^s \oplus \boldsymbol{E}^u$ and $m \in \N$ such that for every $x\in X$
\[  \|(F_m(x))\big|_{\boldsymbol{E}^s(x)}\|<\frac{1}{2}, \quad
     \|(F_{-m}(x))\big|_{\boldsymbol{E}^u(x)}\|<\frac{1}{2}.\]
     
Note that every uniformly hyperbolic cocycle with non-trivial splitting admits some dominated splitting. Moreover, these two notions are equivalent in two-dimensional volume-preserving linear cocycles. 

Recall that we denote the space of cocycles admitting a dominated splitting by $\cD$. When we say a cocycle \emph{stably does not admit dominated splitting}, we mean it is contained in the $\cC^0$-interior of $\cD^c$.
}

\subsection{Symbolic Dynamics}

Given $m\in \bN$, let $\cA_m:=\{1,2,\ldots,m\}$. We denote the set of all bi-infinite sequences in $\cA_m$ by $\Sigma_m$. When $m\in \N$ is fixed, we drop the subscript $m$ in the notations. 
For integers $n<n'$,  and $\ua=(a_{n},a_{n+1},\ldots,a_{n'})\in \cA^{n'-n+1}$, the cylinder $\cylinder[n;\ua]\subset \Sigma_m$ is defined by 
\[\cylinder[n;\ua]:=\{(w_i)_{i\in \Z}\in \Sigma_m  :  w_i=a_i ~ \text{for}~ i=n,\ldots,n'\}.\]
 When $\ua\in \cA^{2n+1}$ for some $n\in \N\cup \{0\}$, we write $\cylinder[\ua]$ instead of $\cylinder[-n;\ua]$. 
In particular when $a\in \cA$, $\cylinder[a]=\cylinder[0;a]$.

 We consider the topology on $\Sigma_m$ generated by the cylinders.
 Moreover, we consider the following metric on  $\Sigma_m$, which induces the mentioned topology. 
\[d\big((x_i)_{i\in \Z},(y_i)_{i\in \Z}\big)=2^{-n},\] where $n= \min \{|i|: x_i \neq y_i \}$. 
The shift map $\sigma:\Sigma_m\to \Sigma_m$ defines a homeomorphism  by  
\[\sigma(\ldots,x_{-1};x_0,x_1,\ldots):=(\ldots,x_0;x_1,x_2,\ldots).\]
 For every $x=(x_i)_{i\in \Z}\in \Sigma_m$, the local stable and unstable set of $x$ is defined by 
 \begin{align*}
     W_{loc}^s(x)&:=\Big\{(y_i)_{i\in \Z} :  y_i=x_i~ \text{for}~ i\geq 0\Big\},\\
       W_{loc}^u(x)&:=\Big\{(y_i)_{i\in \Z} :  y_i=x_i~ \text{for}~ i\leq  0\Big\}.
 \end{align*}

 For $k\in \N$ and $\underline{x}=(x_1,\ldots,x_k)\in \cA^k$, one can associate a unique periodic point to $\underline{x}$, namely $(\ldots, \underline{x};\underline{x},\underline{x},\underline{x},\ldots)$. Note that for every $n\geq 1$, the periodic point associated to $\underbrace{\underline{x}\underline{x}\cdots\underline{x} \underline{x}}_n$ and to  $\underline{x}$ are the same.

An important class of subsystems of $\sigma:\Sigma_m\to \Sigma_m$ are subshifts of finite type. 

\begin{definition}
    For an $m\times m$ matrix  $H$ with $0$ or $1$ entries,   \[\Sigma_{H}:=\big\{(x_i)_{i\in \Z} \in \Sigma_m : \forall i\in \Z~~~ H_{x_i x_{i+1}}=1 \big\}\]
     is a closed subset of  $\Sigma_m$ invariant under $\sigma$. The pair $(\Sigma_H,\sigma) $ is called a \textit{subshift of finite type} or \emph{SFT} for short. 
\end{definition}

\subsection{Linear cocycle on trivial bundle} 

Given a smooth action of a topological group $G$ on $\mathbb{R}^d$, a $G$ cocycle over homeomorphism $f: X \to X$ is a linear cocycle defined on the trivial vector bundle $X \times \mathbb{R}^d$ over $X$ by $F(x,v) = (f(x), A(x)v)$, where $A: X \to G$ is a continuous map. 
For simplicity, we denote such a cocycle by $(f,A)$ instead of the full notation $(X \times \mathbb{R}^d, X, F, f)$. The main class of cocycles studied in this paper are those for which $G$ is a subgroup of $\gldr$ (especially $G=\sldr$ or $G=\gldr$) and the natural action of $d \times d$ matrices on $\R^d$ is assumed. 
In such cases, for $n\in \Z$, the $n$-the iteration of homeomorphism $F$ is a cocycle represented by $(f^n,A_n)$, where for $n\in \N$, 
\begin{align*}
A_n(x)&:=A(f^{n-1}(x))\cdots A(f(x))A(x),  \\
A_{-n}(x)&:=A(f^{-n}(x))^{-1}\cdots A(f^{-1}(x))^{-1},
\end{align*}
and for $n=0$, $A_0$ is the map assigning the identity matrix to each $x\in X$.

\subsection{Hölder regularity}\label{sec:holder}
     A continuous cocycle $(f,A)$ is an $\alpha$-{Hölder} cocycle and  is denoted by $(f,A) \in \mathcal{C}^{\alpha}$ if there exist $C>0$ such that
     $$ \|A(x)- A(y)\| \leq C d(x,y)^\alpha$$
     for every points $x,y \in X$. For every $\eta>0$,the $\cC^\alpha_{loc}$ semi-norm and $\cC^\alpha_{loc}$ norm of an $\alpha$-Hölder cocycle $A$ are defined, respectively by
\begin{align*}
    |A|_{\mathcal{C}^{\alpha},\eta}&:=\sup\limits_{0<d(x,y)<\eta} \frac{\|A(x)-A(y)\|}{d(x,y)^\alpha},\\
     \|A\|_{\mathcal{C}^{\alpha},\eta}&:= \|A\|_{\mathcal{C}^0}+|A|_{\cC^\alpha,\eta}.
\end{align*}
In addition $\cC^\alpha$ semi-norm and $\cC^\alpha$ norm of an $\alpha$-Hölder cocycle $A$ are defined as follows:
\begin{align*}
    |A|_{\mathcal{C}^{\alpha}}&:=|A|_{\mathcal{C}^{\alpha},\mathrm{diam}(X)},\\
     \|A\|_{\mathcal{C}^{\alpha}}&:= \|A\|_{\mathcal{C}^0}+|A|_{\cC^\alpha}.
\end{align*}

One can verify that for every pair of $\alpha$-Hölder cocycles $(f,A),(f,B)\in \cC^\alpha$, \[|A+B|_{\cC^\alpha}\leq |A|_{\cC^\alpha}+|B|_{\cC^\alpha}.\] Therefore, the norm $\|A\|_{\mathcal{C}^{\alpha}}$ can be  used to define a metric on ${C}^{\alpha}$, 
\[d_{\mathcal{C}^{\alpha}}(A,B):=
         \|A-B\|_{\mathcal{C}^{\alpha}}.\]

{Let $\cS$ be subset of $\alpha$-Hölder cocycles on trivial bundle, over $f$. The closure of $\cS$ in the $\cC^\alpha$ topology is denoted by $\mathrm{cl}_{\cC^\alpha}(\cS)$. }

         A special family of Hölder continuous cocycles is locally constant cocycles defined below.

\subsection{Locally constant cocycles}
\label{subsec:loc-const-transition}

{For $f:X\to X$ and $A:X\to \mathrm{GL}(d,\R)$, we say the linear cocycle $(f,A)$ is {\it locally constant}, if there exist $n\in \N$ and disjoint open sets $U_1,\ldots, U_n\subset X$ such that $\bigcup_{i=1}^nU_i=X$ and $A|_{U_i}$ is constant for every $i=1,\ldots,n$.

We remark that under the assumptions above, since    $U_i=X\setminus(\bigcup_{j\neq i} U_j)$,  every $U_i$ is both open and closed. Therefore, for every positive real number  $\eta$  smaller than the distance between $U_i$'s,  one has $A(x)=A(y)$ whenever $x,y\in X$ with $d(x,y)<\eta$.

  \begin{lemma}\label{lem:Hölder}
     Given locally constant cocycle $(f,A)$, let  $\eta>0$  such that $A(x)=A(y)$ for every $x,y\in X$ with $d(x,y)<\eta$. Then, for every $\epsilon >0$, if $d_{\mathcal{C}^{\alpha}}(A,B)<\epsilon$ for some linear cocycle $(f,B)$, then $|B|_{\mathcal{C}^{\alpha},\eta}<\epsilon$, for every $\alpha>0$. 
 \end{lemma}
 \begin{proof}
Since $A(x)=A(y)$ for every $x,y\in X$ with $d(x,y)<\eta$,  we get $|A|_{\cC^\alpha,\eta}=0$ for every $\alpha>0$. If $d_{\mathcal{C}^{\alpha}}(A,B)<\epsilon$, we have $d_{\mathcal{C}^{\alpha},\eta}(A,B)<\epsilon$. Therefore, 
$$|B|_{\mathcal{C}^{\alpha},\eta}\leq |A-B|_{\mathcal{C}^{\alpha},\eta}+|A|_{\mathcal{C}^{\alpha},\eta}\leq \|A-B\|_{\mathcal{C}^{\alpha},\eta}<\epsilon.$$
This finishes the proof of Lemma. 
 \end{proof}
 \begin{corollary}\label{cor:Hölder}
     For a constant cocycle $(f,A)$ and $\epsilon >0$, if $d_{\mathcal{C}^{\alpha}}(A,B)<\epsilon$ for some linear cocycle $(f,B)$, then $|B|_{\mathcal{C}^{\alpha}}<\epsilon$, for every $\alpha>0$.
 \end{corollary}

Throughout this paper, many statements are formulated for the locally constant cocycles over the shift map $\sigma:\Sigma_m\to \Sigma_m$.

 {Given  integers $n,n'\geq 0$ and  $\phi:\cA^{n+n'+1}\to \mathrm{GL}(d,\R)$, we say that a locally constant $\mathrm{GL}(d,\R)$ cocycle $(\sigma,A)$ over  shift is {\it represented by $\phi$}, if for every $\ua\in \cA^{n+n'+1}$, 
 \[A|_{\cylinder[-n;\ua]}=\phi(\ua)\in \mathrm{GL}(d,\R).\]
 Moreover, we say that a locally constant cocycle $(\sigma,A)$ depends on positions in ${\ldbrack -n,n'\rdbrack}:=\{-n,\ldots,n'\}$, if there exists  $\phi:\cA^{n+n'+1}\to \mathrm{GL}(d,\R)$ so that  $A$ is represented by $\phi$. 

 When $n=n'=0$, we say that the locally constant cocycle depends only on the {\it {zeroth position}}. {A correspondence can be established between locally constant cocycles depending on zeroth position and the iterated function systems (IFS) generated by the matrices in the range of $\phi$. By an IFS  generated by $\langle \phi_1,\ldots,\phi_m \rangle^+$, we mean a locally constant cocycle depending on the zeroth position over the shift map $\sigma:\Sigma_m \to \Sigma_m$. }

The proposition below provides a straightforward characterization of locally constant cocycles over the shift, asserting that every locally constant cocycle can be represented in a specified form, offering an equivalent definition for locally constant cocycles over the shift.

 \begin{proposition} \label{prop:lc-over-shift}
     For every locally constant cocycle  $(\sigma,A)$ there exists $k\in \N$ such that $A$ depends on positions in ${\ldbrack-k,+k\rdbrack}$. 
 \end{proposition}
 }

    Denote the set of all finite words of elements of $\cA$ by $\cA^{\mathrm{fin}}$. Clearly $\cA^{\mathrm{fin}}=\cup_{n\in \N}\cA^{n}$. For every $\uw\in \cA^{\mathrm{fin}}$, we say that $\uw$ has length $n$ and write $|\uw|=n$, if $\uw\in \cA^n$. 

We consider two operations on the set $\cA^{\mathrm{fin}}$. For $w=(w_1,\ldots,w_n)\in \cA^n$ and $w'=(w'_1,\ldots,w'_{n'})\in \cA^{n'}$, we define the concatenation 
\begin{equation}\label{eq:def-concatenation}
    \uw\,\uw':=(w_1,\ldots,w_n,w'_1,\ldots,w'_{n'})\in \cA^{n+n'}.
\end{equation}

Also, for every positive integer $t\leq n,n'$, if $(w_{n-t+1},\ldots,w_n)=(w'_1,\ldots,w'_t)$, then we define $ \underline{w}\vee_{t} \underline{w'}\in \cA^{n+n'-t}$ by 
\begin{equation}\label{eq:def-merge}
  \underline{w}\vee_{t} \underline{w'}:=(w_1,\ldots,w_n,w'_{t+1},w'_{t+2},\ldots,w'_{n'})=(w_1,\ldots,w_{n-t},w'_1,\ldots,w'_{n'}).
\end{equation}
If there is no ambiguity, we omit dependence to $t$ and write $\uw\vee \uw'$ instead of $\uw\vee_{t} \uw'$.  Note that for the special case $t=0$, $\uw\vee_0 \uw'=\uw\,\uw'$.

 \begin{definition}\label{def:transition-locall}
Given $l,l'\in \N$ and $\ua,\ub\in \cA^{\mathrm{fin}}$ with $|\ua|=l$, $|\ua'|=l'$, by  a {\it transition}  from $\ua$ to $\ua'$, we mean a word 
 $\underline{w}=(w_1,w_2,\ldots,w_n) \in \cA^{\mathrm{fin}}$ (with $|\uw|\geq l,l'$), such that 
  \[\ua=(w_1,\ldots,w_{l}) , \quad \ub=(w_{n-l'+1},\ldots,w_{n}). \] 
 The set of all transitions from $\ua$ to $\ua'$ is denoted by $T(\ua,\ua')$. 
 \end{definition}

\begin{remark}
    One can give the following remarks about the definition of transition. 
     \begin{itemize}
  \item If $\uw\in T(\ua,\ua')$ satisfies $|\underline{w}|>|\ua|+|\ua'|$, then there exists $\ub\in \cA^{\mathrm{fin}}$ such that $\uw=\ua\,\ub\,\ua'$.
 \item If $\ua,\ua',\ua''\in \cA^{\mathrm{fin}}$ and  $\uw\in T(\ua,\ua'),\uw'\in T(\ua',\ua'')$, then $\uw\vee_{|\ua'|}\uw'\in T(\ua,\ua'')$.
  \end{itemize}
  \end{remark}

An equivalent way of looking at transitions is to define them as walks\footnote{A walk is a sequence of edges directed in the same direction which joins a sequence of vertices.} in certain graphs associated with SFTs. For  a locally constant cocycle $A$ over a subshift,  depending on positions in $\ik$, 
One can associate a directed graph {(the so-called \emph{De Bruijn graph})} whose vertices correspond to elements of  $\cA^{2k+1}$ and for $\ua,\ub\in \cA^{2k+1}$ there is a directed edge from $\ua$ to $\ub$ if and only if $b_i=a_{i+1}$ for every $i\in {\ldbrack -k,k-1\rdbrack}$ (see Figure \ref{fig:graph}). In this viewpoint, a transition from $\ua\in \cA^{\ik}$ to $\ub\in \cA^{\ik}$ is a walk in this graph starting and ending in the vertices corresponding to $\ua$ and $\ub$, respectively. 

\begin{figure}[t]
     \centering
     \begin{subfigure}[b]{0.3\textwidth}
         \centering         \includegraphics[width=\textwidth]{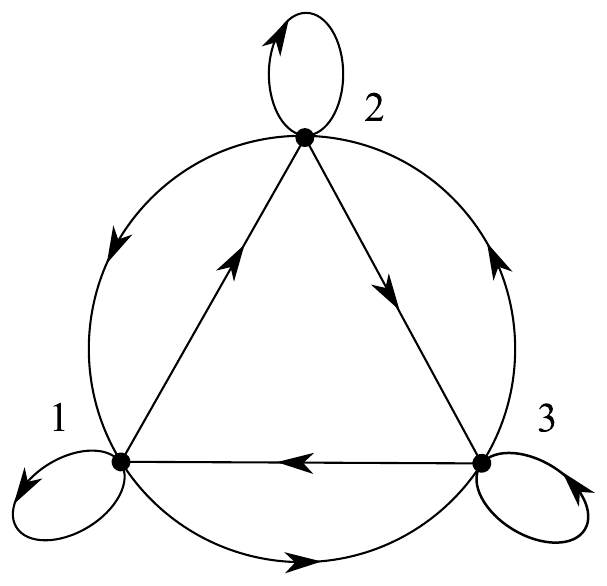}
     \end{subfigure}
     \qquad\qquad
     \begin{subfigure}[b]{0.4\textwidth}
         \centering
\includegraphics[width=\textwidth]{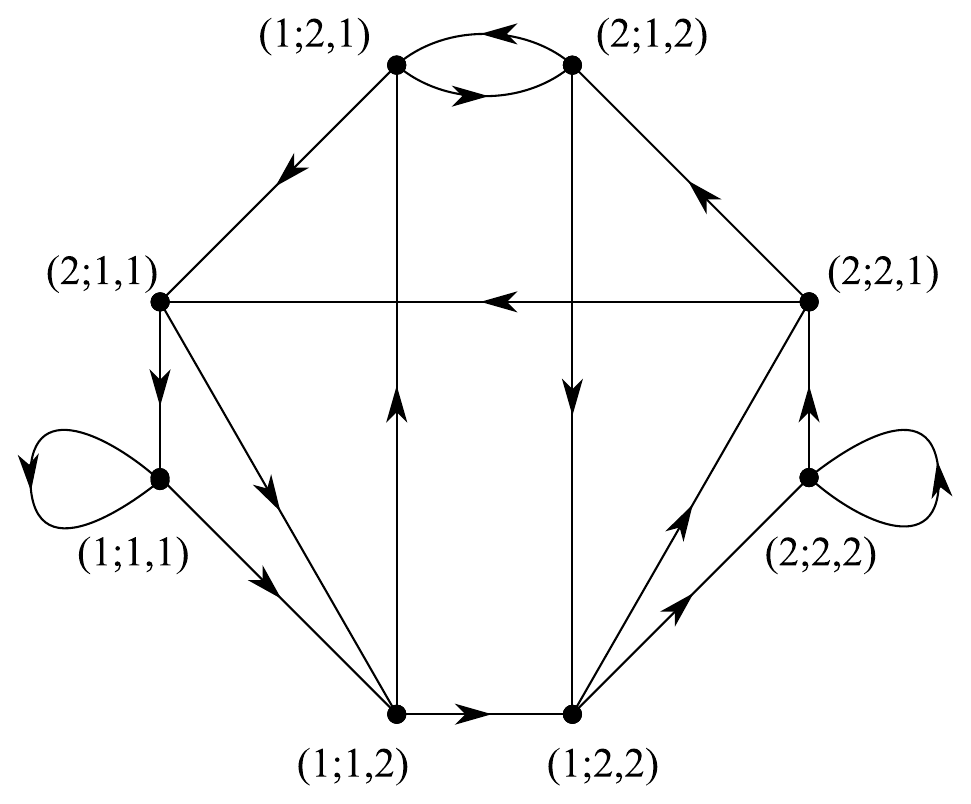}
     \end{subfigure}
    \caption{The left and right figures are the transition graphs for the locally constant cocycles ${A}:\Sigma_3\to \sldr$ and  $B:\Sigma_2\to \sldr$  depending on zeroth position, and positions in $\ldbrack -1,1\rdbrack$, respectively.}
    \label{fig:graph}
\end{figure}

We can also define matrix transition in line with the definition of transition.
\begin{definition}
  Let $[T']$ and $[T'']$ be two finite words of matrices of length $\ell'$ and  $\ell''$, respectively. By a matrix transition from $[T']$ to $[T'']$, we mean a finite word $[T]=(T_1,\ldots,T_\ell)$ of matrices, with $\ell \geq \ell', \ell''$, such that
    \[[T']=(T_1,\ldots,T_{\ell'}) , \quad [T'']=(T_{\ell-\ell''+1},\ldots,T_{\ell}).\] 
\end{definition}

 Now, let $(\sigma,A)$ be a locally constant $\mathrm{GL}(d,\R)$ cocycle over the shift map, depending on positions in $\ik$ and  represented by $\phi:\cA^{2k+1}\to \mathrm{GL}(d,\R)$.
 Then, for every transition $\underline{w}=(w_1,w_2,\ldots,w_n)$ from $\ua\in \cA^{2k+1}$ to $\ua'\in \cA^{2k+1}$, it follows that the word $\big(\phi(\underline{w}_{n-2k-1}),\dots ,\phi(\underline{w}_2),\phi(\underline{w}_1)\big)$ is a matrix transition from $(\phi(\underline{a}))$ to $(\phi(\underline{a'}))$,  where $\underline{w}_i=(w_i,w_{i+1},\ldots,w_{i+{2k}})$. Note that here, $(\phi(\underline{a}))$ and $(\phi(\underline{a'}))$ are two words of length one. In this case, we denote 
 \[A_{\transition{\uw;2k+1}}:=\phi(\underline{w}_{n-2k-1})\cdots \phi(\underline{w}_2)\phi(\underline{w}_1).\]
 Note that $2k+1$ is the length of words $\uw_i$. So whenever the length $2k+1$ is clear from the context, we omit it from the notation and write $A_{\transition{\uw}}$ instead of $A_{\transition{\uw;2k+1}}$.  
Notice that if $\uw\in T(\ua,\ua')$ and  $\uw'\in T(\ua',\ua'')$, then
\[A_{[\uw\vee_{2k+1} \uw'\rangle}=A_{\transition{\uw'}}A_{\transition{\uw}}.\]

\begin{lemma}\label{lem: locall-realization}
 For a locally constant $\mathrm{GL}(d,\R)$ cocycle $(\sigma,A)$ over the shift map, depending on positions in $\ik$, let $p',p''$ be two periodic points associated to $\ua',\ua''\in \cA^{2k+1}$, respectively. Then, for  every transition $\uw\in T(\ua',\ua'')$, there exists a periodic point $p$ such that
 \[A_{|\uw|-2k-1}(p)=A_{\transition{\uw }}.\]
    
\end{lemma}

For the proof, it is enough to consider $p$ to be the periodic point such that $\sigma^k(p)$ is associated with a word in $\cA^{\rm fin}$ starting with $\uw$.

\section{Cocycles with quasi-conformal orbit}\label{sec:qc and bounded}
{
Recall that for the homeomorphism $f:X\to X$ of the compact metric space $X$ and  continuous map $A:X\to \gldr$, we say that the cocycle $(f,A)$ admits quasi-conformal orbits if there exists $x\in X$ such that {\emph{the condition number} (i.e. ratio of norm and conorm) of all matrices $A_n(x)$ ($n\in \Z$) is bounded by a constant independent of $n$, in other words,}
\begin{equation}\label{def:qc-orb-Z}
    \sup_{n\in \Z} \|A_n(x)\|.\|A_n(x)^{-1}\|<+\infty.
\end{equation}
For $\sldr$ cocycles, the above condition is equivalent to the existence of $x\in X$ satisfying 
\begin{equation*}
    \sup_{n\in \Z} \|A_n(x)\|<+\infty.
\end{equation*}
Moreover, we say $(f,A)$ admits fiberwise
bounded orbits if there exists $x\in X$ such that 
\begin{equation}\label{def:bd-orb-Z}
    \sup_{n\in \Z}\|A_n(x)^{\pm 1}\|<+\infty.
\end{equation}
We denote the set of all continuous cocycles $(f,A)$ admitting quasi-conformal orbits by $\qc$,  and the set of those cocycles admitting fiberwise
bounded orbits by $\cB$. It is easy to verify that $\cB\subset \qc$, and if $(f,A)$ is a $\sldr$ cocycle,  $(f,A)\in \cB$ if and only if $(f,A)\in \qc$. }

\begin{remark}
    We can also define cocycles admitting quasi-conformal or fiberwise bounded forward orbits (resp. backward orbits), by replacing $n\in \Z$ in \eqref{def:qc-orb-Z} and \eqref{def:bd-orb-Z} by $n \in \N$ (resp. $-n\in \N$). 
\end{remark}

\subsection{Quasi-conformal orbits of a cocycle.}

Throughout this section, $X$ is a compact metric space and $f:X\to X$ is a homeomorphism, and $A:X\to \gldr$ is a continuous map. All the results are formulated for the 
 cocycles $(f,A)$ over trivial bundles.

For every $\kappa\geq 1$, we denote 
\[{\rm QC}^+(A,\kappa):=\big\{x\in X\,:\, \|A_n(x)\|.\|A_n(x)^{-1}\|\leq \kappa~ \text{for all }n \in \N\big\}.\]

Clearly, the set of points $x\in X$ admitting quasi-conformal (forward) orbits for $A$, coincides the union $\bigcup_{\kappa\geq 1} {\rm QC}^+(A,\kappa)$.
The next lemma provides the basic properties of sets ${\rm QC}^+(A,\kappa)$.  
{\begin{lemma}\label{lem:qc-properties}
     For every  $\kappa\geq 1$,
     \begin{itemize}
         \item[(i)] If for some non-negative integers $m<n$, 
         $\|A_m(x)\|.\|A_m(x)^{-1}\|\leq \kappa$ and 
         $\|A_n(x)\|.\|A_n(x)^{-1}\|\leq \kappa$, then 
         $\|A_{n-m}(f^m(x))\|.\|A_{n-m}(f^m(x))^{-1}\|\leq \kappa^2$.
         \item[(ii)] ${\rm QC}^+(A,\kappa)$ is closed. 
         
    \item[(iii)]  If $x\in {\rm QC}^+(A,\kappa)$, then $\cO^+(x)\subset {\rm QC}^+(A,\kappa^2)$. 
    \item[(iv)] If $x\in {\rm QC}^+(A,\kappa)$, then for every $y\in \omega(x)$, and every $n\in \Z$
\[\|A_n(y)\|.\|A_n(y)^{-1}\|\leq \kappa^2.\]
In particular,  $y\in {\rm QC}^+(A,\kappa^2)$.
     \end{itemize}
\end{lemma}
\begin{proof}
    Part (i): Since $A_{n-m}(f^m(x))=A_n(x)A_m(x)^{-1}$, 
    \[\|A_{n-m}(f^m(x))\|.\|A_{n-m}(f^m(x))^{-1}\|\leq  \|A_n(x)\|.\|A_m(x)^{-1}\|.\|A_n(x)^{-1}\|.\|A_m(x)\|\leq \kappa^2.\]
    Part (ii) is a consequence of continuity of $A_n$; and 
    Part (iii) follows immediately from part (i). 
   
To prove part (iv),  fix an element of $\omega(x)$, say $y$ and $n\in \Z$.  Then, there exists a sequence $n_1<n_2<\cdots$ of positive integers such that 
    \begin{equation}\label{eq:omega-limit}
       \lim\limits_{j\to\infty}f^{n_j}(x)=y.
    \end{equation}
     Now, for sufficiently large $j$, $n_j+n>0$ and for such $j$'s, using Part (i), we get 
     \(\|A_n(f^{n_j}(x))\|.\|A_n(f^{n_j}(x))^{-1}\|\leq \kappa^2\). 
     In view of \eqref{eq:omega-limit}, and since $A_n$ and $A_n^{-1}$ are both continuous, we conclude that $\|A_n(y)\|.\|A_n(y)^{-1}\|\leq \kappa^2$. 
\end{proof}
}

The statement of Lemma \ref{lem:qc-properties} holds for the set of points whose fiberwise forward orbit is bounded, with a similar proof. More precisely, 

\begin{lemma}\label{lem:bounded-properties}
     If for $\kappa\geq 0$, we denote 
\[{\rm B}^+(A,\kappa):=\{x\in X\,:\, \|A_n(x)^{\pm 1}\|\leq \kappa~ \text{for all }n \in \N\}.\]
Then,
   \begin{itemize}
       \item[(i)] If for some non-negative integers $m<n$, $\|A_m(x)^{\pm 1}\|\leq \kappa$ and $\|A_n(x)^{\pm 1}\|\leq \kappa$, then $\|A_{n-m}(f^m(x))\|\leq \kappa^2$. 
    \item[(ii)]  ${\rm B}^+(A,\kappa)$ is closed. 
    \item[(iii)]  If $x\in {\rm B}^+(A,\kappa)$, then $\cO^+(x)\subset {\rm B}^+(A,\kappa^2)$.  
\item[(iv)]  If $x\in {\rm B}^+(A,\kappa)$, then for every $y\in \omega(x)$ and every $n \in \Z$
\[\|A_n(y)^{\pm 1}\| \leq \kappa^2.\]
In particular,  $y\in {\rm B}^+(A,\kappa^2)$.
   \end{itemize}
\end{lemma}

\begin{remark}
    The statement of Lemma \ref{lem:qc-properties} and \ref{lem:bounded-properties}, with the same proof, holds if one replaces the set ${\rm QC}^+(A,\kappa)$ and ${\rm B}^+(A,\kappa)$ with the set of points $x\in X$ for which the norm of cocycle is quasi-conformal or bounded by $\kappa$ for every point in the full orbit of $x$, respectively.
\end{remark}

\subsection{Non-stability of quasi-conformal orbits}
In this subsection, we consider linear cocycles on trivial bundles over an arbitrary base map (not necessarily the shift map). We will show that one can violate the quasi-conformality of a specific orbit by a small $\cC^0$ perturbation of the cocycle.

{
\begin{lemma}\label{lem:stable qc fixed orbit}
For every homeomorphism $f:X\to X$, and every $x_0\in X$, 
\begin{itemize}
  \item If $d> 2$ and $x_0$ is periodic for $f$, then the set of $\sldr$ cocycles $(f,B)$ for which  $\{\|B_n(x_0)\|\}_{n\in \N}$  is unbounded forms an open and dense subset in the $\mathcal{C}^0$ topology. 
\item If $d\geq  2$ and $x_0$ is not periodic for $f$, then the set of $\sldr$ cocycles $(f,B)$ for which  $\{\|B_n(x_0)\|\}_{n\in \N}$ is unbounded forms a residual (hence, dense) subset in the $\cC^0$ topology. 
\end{itemize}
\end{lemma}
}

\begin{proof}
If $d>2$ and $x_0$ is a periodic point, then the conclusion follows from Proposition \ref{prop:generic hyp periodic} applied to $X=\cO(x_0)$, and Lemma \ref{lem:extension}. 

Now, let $x_0$ be a non-periodic point and let $d\geq 2$.
For every $k\geq 1$, define \[\cB_\kappa:= \big\{B:X \to \sldr:  \sup\limits_{n\in \N} \|B_n(x_0)\| \leq \kappa\big\}.\]
It is clear that for every $k\geq 1$, $\cB_\kappa$ is closed. On the other hand, the set of  continuous maps $B:X\to \sldr$, for which $\{\|B_n(x_0)\|\}_{n\in \N}$ is unbounded, equals to  \(\bigcap_{\kappa \in \N} \cB_\kappa^c\).

So, it suffices to prove that for every $\kappa\geq 1$, the $\cC^0$-interior of $\cB_\kappa$ is empty. This implies the lemma as the set of $\sldr$ cocycles over $f$ is a Baire space. 

\medskip

   \noindent
    \textbf{Claim.} The interior of $\cB_\kappa$ is empty.

     Let $\epsilon_0 >0$ and $B \in \cB_\kappa$. For every  $r>1$, let $D_r:=\mathrm{Diag}(r,r^{-\frac{1}{d-1}},\ldots,r^{-\frac{1}{d-1}})$ be a diagonal matrix and for every vector $v \in \R^d$ choose $R_v \in \mathrm{SO}(d,\R)$ such that $R_v(e_1)={v}/{|v|}$, where $(e_1,\ldots,e_d)$ is the standard basis for $\R^d$. 
     
     Consider an arbitrary unit vector $v_0 \in \R^d$ and $n \in \N$. Now, for every $i=1,2,\ldots,n$, 
     denote $v_i:= B_i(x_0)(B_n(x_0))^{-1}v_0$ and  define \[\hat{B}(f^i(x_0))=R_{v_i}D_rR^{-1}_{v_i}B(f^i(x_0)).\]
     Notice that as $x_0$ is not periodic, $\hat{B}$ is well-defined (for every $n$). 
     Since $v_n=v_0$, then $\hat{B}_n(x_0)(B_n(x_0))^{-1}v_0=r^nv_0$. Consequently,
     \begin{align}\label{inequality}
     \| \hat{B}_n(x_0)\|  \geq \frac{r^n}{\|B_n(x_0)^{-1}\|}\geq \frac{r^n}{\|B_n(x_0)\|^{d-1}}\geq r^n \kappa_0^{-(d-1)}. 
     \end{align}
     To complete the proof of claim,
     let $\delta$ be the positive number given by Lemma \ref{lem:extension} so that every $\delta$-perturbation of cocycle $B$ on a compact set can be extended to a cocycle on $X$ $\epsilon_0$-close to $B$.

     Now, take $r>1$ close enough to $1$, so that $\|D_r-\id\|.\|B\|<\delta$. Moreover, we can take $n$ large enough to have $r^n \kappa_0^{-(d-1)} > k$. We take $K$ in Lemma \ref{lem:extension} equal to $\{x_0,\ldots, f^{n-1}(x_0)\}$. 
     According to Lemma \ref{lem:extension},  there exist $\Tilde{B}:X\to \sldr$, with  $d_{\cC^0}(\Tilde{B},B)<\epsilon_0$ and $\Tilde{B}(f^i(x_0))=\hat{B}(f^i(x_0))$. Moreover, \eqref{inequality} implies that  $\Tilde{B} \notin \cB_\kappa$ and the proof of claim and lemma is complete. 
    \end{proof}
    }

    {
    Using Lemma \ref{lem:stable qc fixed orbit},
    we get the following
    corollary regarding the interior of the set of cocycles admitting quasi-conformal orbits, when  the base system is minimal. }

\begin{corollary}\label{cor:minimal}
{
For every $d>2$ and every minimal homeomorphism $f:X\to X$, the $\cC^0$-interior of the set of $\sldr$ cocycles admitting quasi-conformal orbits is empty.}
\end{corollary}

   \begin{proof}
   Observe that for every $A:X \to \sldr$, if  ${\rm QC}(A,\kappa)\neq \emptyset$ for some  $\kappa\geq 1$, then, using Lemma \ref{lem:qc-properties},  and minimality of $f$,  ${\rm QC}(A,\kappa^2)=X$.

  This implies that if there exists a $\cC^0$-open set of $\sldr$ cocycles over $f$ which admit a quasi-conformal orbit, then for every $x\in X$, the set of cocycles with bounded forward orbit at $x$, has non-empty $\cC^0$-interior. But due to  Lemma  \ref{lem:stable qc fixed orbit}, this is impossible.
\end{proof}

{The following lemma clarifies the sensitivity of $\gldr$ cocycles with bounded orbits to the natural perturbations. {The proof is straightforward and is left to the reader. 
} 
\begin{lemma}
      Let $(f,A)$ be a $\gldr$ cocycle and $Y\subset X$ be an $f$-invariant compact set. Assume that for every $y\in Y$, the sequence $\{\|A_n(y)\|^{\pm 1}\}_{n\in \Z}$ is bounded.
      If $\lambda\in \R\setminus \{0,\pm 1\}$, then the cocycle  $( f\big|_Y,\lambda A)$ does not admit fiberwise bounded orbits.
\end{lemma}

}

\section{Covering condition for cocycles}\label{Sec:covering}

In this section, we introduce a covering criterion that stably violates dominated splitting, based on the recent work \cite{FNR} on the stable ergodicity of smooth expanding actions. The covering condition in \cite{FNR} ensures the existence of quasi-conformal orbits for (pseudo)-semigroup action.

\subsection{Covering condition in linear groups}\label{subsec:covering-linear-groups}
In this section, we recall and adopt the covering condition in \cite{FNR} for the action of subgroups of  $\gldr$ on themselves and explore its implications for identifying quasi-conformal sequences of matrices.

    For $\kappa\geq 1$, we say a sequence $\{D_i\}_{i=1}^\infty$ in $\gldr$   is \emph{$\kappa$-conformal} if and only if for every $n\in \N$
    \[ \|D_n\cdots D_1\|.\|(D_n\cdots D_1)^{-1}\|\leq \kappa.\]
    Also, the sequence is \emph{quasi-conformal} if it is $\kappa$-conformal for some $\kappa\geq 1$. 

For $\kappa\geq 1$ and $\cG=\{D_1,D_2,\ldots,D_m\}\subset \gldr$, we say $\langle \cG \rangle^+$, the semigroup generated by $\cG$, has a \emph{$\kappa$-conformal branch}, if there exists a $\kappa$-conformal sequence in $\cG$. 

    We need the following Lemmas from \cite{FNR}. Lemma \ref{lem:U-branch} provides an equivalent condition for the existence of a quasi-conformal sequence from a given set. 

    \begin{lemma}{\cite[Lemma 3.5]{FNR}}\label{lem:U-branch}
    Let $\cG\subseteq \sldr$ be a finite set. Then, there exists $\kappa>1$ such that $\langle{\cG}\rangle^+$ has a  $\kappa$-conformal branch if and only if $ ~\overline{\cU}\subseteq \cG^{-1}\cU$ for some $\cU\subseteq \sldr$ with compact closure. Moreover, if $\cU$ is open,  $\langle{\cG}\rangle^+$ robustly has a $\kappa$-conformal branch for some $\kappa>1$.   
    \end{lemma}

The lemma, in particular, implies that if all the elements of $\cG$ have a common expanding eigenvector, or more generally, if the locally constant cocycle associated with $\cG$ admits a dominated splitting, then $\cG$ does not satisfy the covering condition (for any open set with compact closure). Indeed,  dominated splitting or common expanding eigenvector implies that every branch in $\cG$ is unbounded. See \S \ref{subsec:question-IFS}  for further discussion. 

In view of the statement of the above lemma, for topological group $G$, we say that the finite set $\cG=\{g_1,\ldots,g_m\}\subset G$ satisfies the \emph{covering condition} if there exists a non-empty open set $\cU\subset G$ with compact closure such that 
\begin{equation}\label{eq:covering-IFS}
    \overline{\cU}\subset \cG^{-1}\cU=\bigcup\limits_{i=1}^mg_i^{-1}\cU.
\end{equation}
In the next subsection, in Definition \ref{def:gen-cov}, the covering condition will be generalized to the locally constant cocycles over the shift.

The next lemma guarantees the existence of finite subsets of $\sldr$ with exactly $d^2$ elements and with elements sufficiently close to $\id$, satisfying the covering condition.

\begin{lemma}{\cite[Lemma 3.8]{FNR}}\label{lem:sldr-covering}
    For any neighborhood $\cU_0$ of the identity in $\sldr$, there exists an open set $\cU\subseteq \cU_0$ and a finite set $\cG\subseteq \cU_0$ with  $d^2$ elements such that $\overline{\cU}\subseteq \cG^{-1}\cU$.
    \end{lemma}
    {We can also satisfy the covering condition in $\gldr$. let $\{D_1,\ldots,D_{d^2}\}\subset\sldr$, satisfy the covering condition for some open set $\cU\subset\sldr$. Consider two maps $t\mapsto 2t$ and $t\mapsto t/3$ on $\R$.  
It is easy to see that these two maps satisfy the covering condition for any bounded interval $(a, b)\subset \R$ provided that  $b>6a$. So, we can conclude that the family 
$\{2D_1,\ldots,2D_{d^2},\frac{1}{3}D_1,\ldots,\frac{1}{3}D_{d^2} \}$, satisfy the covering condition for the open set $\cU':=\{\lambda D : a<\lambda^d<b\And A \in \cU \}$ in $\gldr$. In general, we can apply the following corollary of Lemma \ref{lem:sldr-covering} to satisfy the covering condition in $\gldr$. 

\begin{corollary}\label{cor:gldr-covering}
    For any neighborhood $\cU_0$ of the identity in $\gldr$, there exists an open set $\cU\subseteq \cU_0$ and a finite set $\cG\subseteq \cU_0$ with  $2d^2$ elements such that $\overline{\cU}\subseteq \cG^{-1}\cU$.
\end{corollary}
}

The next lemma comes from \cite{FNR} and \cite{homburg_nassiri_2014}.

 \begin{lemma}\label{lem:smaller-U}
 Let $\cG\subset\sldr$ be such that for some non-empty open subset $\cU\subset \sldr$ with compact closure, $\overline{\cU}\subset \cG^{-1}\cU$. Then, there exists $\delta>0$ such that
     $$\overline{\cU}\subset \cG^{-1}\cU_{(\delta)},$$
      where $\cU_{(\delta)}:=\big\{D :D'\in \cU\; \text{for every } D'\in \sldr \text{with }\|D-D'\|<\delta\big\}$. 
 \end{lemma}
 \begin{proof}
   If such a $\delta$ did not exist, then for every positive integer $k$, we would have $\Lambda_k := \overline{\cU} \setminus\big( \cG^{-1}\cU_{(\frac{1}{k})}\big) \neq \emptyset$. Since $\{\Lambda_k\}_{k=1}^\infty$ is a nested sequence of non-empty compact sets, it follows that \[\bigcap\limits_{k=1}^\infty \Lambda_k=\overline{\cU} \setminus\big( \cG^{-1}\cU\big)\neq \emptyset,\]
   which contradicts the assumption $\overline{\cU}\subset \cG^{-1}\cU$.
 \end{proof}

 \begin{lemma}[covering implies immediate covering]
 \label{lem:immediate-covering}
Let $\cH$ be the semigroup generated by $H_1,H_2,\ldots,H_m\in\sldr$.  Assume that there is an open subset $\cU\subset \sldr$ with compact closure such that $\overline{\cU}\subset \bigcup_{H\in \cH} {H}^{-1}\cU$.
Then there exist an open subset $\cV\subset \sldr$ with compact closure such that
$$\overline{\mathcal{V}}\subset \bigcup_{i=1}^m {H_i}^{-1}\mathcal{V}.$$
 \end{lemma}
\begin{proof}
Since $\overline{\cU}$ is compact and it is covered by the union  $\bigcup_{H\in \cH} {H}^{-1}\cU$ of open sets, it can be covered by the union of finitely many of them. So, $\bigcup_{H\in \cF} {H}^{-1}\cU$, for some finite family 
$\cF\subset \cH$.
{We may assume that $\cF$ is the set of all words of length at most $N\in \N$.}
By Lemma \ref{lem:smaller-U} there exists $\delta>0$ such that $\overline{\cU}\subset \bigcup_{H\in \cF} {H}^{-1}\cU_{(\delta)}$.  
Since linear maps are uniformly continuous {on bounded sets}, there exist positive real numbers $\epsilon_0 <\ldots <\epsilon_N<\delta/2$ such that for every 
{$i\in \{1, \dots, m\}$, $j\in \{1, \dots, N\}$, and 
$X, Y \in \bigcup_{H\in \cF} {H}^{-1}\cU$,} 
if $\|X-Y\|< \epsilon_{j-1}$ then $\|H_iX-H_iY\|<\epsilon_j/2$. 
Denote $\epsilon:=\frac{1}{2}\min \{\epsilon_0,\ldots,\epsilon_N\}$.

Fixing $u\in \overline{\cU}$, there exist positive integers  $k\leq N$, and $i_1,\ldots, i_k \leq m$ such that $\mathrm{B}_{\delta}(H_{i_k} \cdots H_{i_1} u) \subset \cU$. Let $u_0:=u$, $\cV_0:=\mathrm{B}_{\epsilon_0}(u_0)$, and $u_j:=H_{i_j}u_{j-1}$,  $\cV_j:=\mathrm{B}_{\epsilon_j}(u_j)$ for $1\leq j \leq k$. 
It follows that  $H_{i_j}\mathrm{B}_{\epsilon}(\cV_{j-1})\subset \cV_j$ and 
$\mathrm{B}_{\delta/2}(\cV_k) \subset \cU$.
Define $\cV_u:= \bigcup_{j=0}^k \cV_j$. 

Hence, the bounded open set $\cV:= \bigcup_{u\in \overline{\cU}}\cV_u$ satisfies the covering property for $H_1,\ldots,H_m$.
\end{proof}

\subsection{Covering condition for cocycles}
\label{subsec:covering cocycles}
{
In this section, we extend the covering condition to include general locally constant cocycles. The concept of a cocycle with covering can be defined for cocycles taking values in more general topological groups $G$. However, for the purposes of this paper, we restrict our attention to cocycles taking values in closed subgroups of linear groups. The most important cases that we will consider in the next sections are when 
$G=\gldr$ or $G=\sldr$.
}

\begin{definition}[Cocycle with covering]\label{def:gen-cov}
Let $d\geq 2$ and  $G$ be a {closed} subgroup of $\gldr$. 
We say that a locally constant $G$ cocycle $(\sigma,A)$ over the shift map, depending on positions in $\ldbrack -k,k\rdbrack$ and represented by $\phi:\cA^{2k+1}\to G$
satisfies the  \emph{covering condition} if there exist
\begin{itemize}
    \item a non-empty open set $\cU\subset G$ with  compact closure,
    \item and a subset $\cA'\subset \cA^{2k+1}$,
\end{itemize}
such that for every $\ua\in \cA'$
\begin{equation}\label{eq:gen-covering}
      \overline{\cU} \subset \bigcup\limits_{\ub \in  \cA'}\bigcup\limits_{\uw\in T(\ua,\ub)}(A_{[\uw\rangle})^{-1}\cU.
\end{equation} 
\end{definition}

\begin{remark} We can give a few remarks on the definition of covering property for cocycles. 
\label{rmk:better-gen-cov}
\begin{itemize}
    \item[1)] If the covering condition \eqref{eq:gen-covering} holds for some open set $\cU$, we may assume that $\cU$ contains the identity matrix, $\id$. Indeed, we can fix an arbitrary element $H_0\in \cU$ and replace $\cU$ in \eqref{eq:gen-covering} with $H_0^{-1}\cU$. 
    \item[2)]  Since $\cU$ is open and $\overline{\cU}$ is compact, according to Lemma \ref{lem:smaller-U}, we deduce that Definition \ref{def:gen-cov} is equivalent to the following: \\
    There exists $\delta>0$ such that for every $\ua\in \cA'$, there is a finite set  $\Lambda_{\ua}$ of transitions of the form $T(\ua,\ub)$ for $\ub\in \cA'$ such that 
    \begin{equation}\label{eq:gen-covering2}
        \overline{\cU} \subset \bigcup\limits_{\uw\in  \Lambda_{\ua}}(A_{[\uw\rangle})^{-1}\cU_{(\delta)}, \quad \text{for every } \ua \in \cA'.
    \end{equation}
\end{itemize}
   
\end{remark}

Note that the Definition \ref{def:gen-cov} extends the covering condition introduced in the previous subsection in \eqref{eq:covering-IFS}. Indeed, a locally constant cocycle $(\sigma,A)$ over $\Sigma_m$ depending on zeroth position can be represented by an $m$-tuple  $(D_1,\ldots,D_m)$ of matrices in $\sldr$. If $(\sigma,A)$ satisfies the covering condition in the sense of Definition \ref{def:gen-cov} for   $\cA'=\cA$ and for some open set $\cU\subset \sldr$ with compact closure, it means that 
\[ \overline{\cU}\subset \bigcup\limits_{D\in \langle \cG\rangle^+} D^{-1} \cU.\]
This because of Lemma \ref{lem:immediate-covering} implies that there exists open set $\cV\subset \sldr$ with compact closure for which we have 
\[ \overline{\cV}\subset \bigcup\limits_{D\in  \cG} D^{-1} \cV,\]
This is exactly the covering condition in \eqref{eq:covering-IFS}. 
{One can construct more sophisticated examples of locally constant cocycles over the shift map which do not depend only on zeroth position, satisfying the covering condition
(see the proof of \cite[Lemma 5.17]{Reshadat-thesis}).
}

\begin{remark}{
Let $(\sigma,A)$ be the locally constant cocycle over the shift map, depending on positions in $\ldbrack -k,k\rdbrack$ and represented by $\phi:\cA^{2k+1}\to \sldr$. If the family $\phi(\cA^{2k+1})$ of matrices satisfies the covering condition with respect to some open set $\cU\subset\sldr$, then it does not necessarily imply \eqref{eq:gen-covering}. Indeed, it is not difficult to construct an example of a uniformly hyperbolic locally constant cocycle over the full shift $\Sigma_2$, represented by some $\phi:\cA^3\to \sldr$ so that the semigroup generated by $\phi(\cA^3)$ is dense in $\sldr$ (see \cite[Example 5.16]{Reshadat-thesis}). }\end{remark}

\section{Criterion for stable quasi-conformality}\label{sec:stability-qc}

{
In this section, we consider locally constant linear cocycles over the shift map and establish the existence of quasi-conformal orbits and their stability in the $\cC^\alpha$ topology under the covering condition introduced in the previous section.

For the case of locally constant cocycles, we first prove the results when the cocycle depends only on the zeroth position and then give the proof for general locally constant cocycle. While the proof for general locally constant cocycles in \S \ref{subsec:stability-qc} gives the results of \S \ref{subsec:stability-qc} as a particular case, we initially present the results for the cocycles depending only on zeroth position to make the idea more clear and improve readability. 
}

 \subsection{Cocycles depending on zeroth position}\label{subsec:stablity-qc-zero-position}
 
As mentioned above, our focus in this subsection is on the perturbation of IFS, namely cocycles that depend on the zeroth position. We will show that for a subgroup $G$ of $\gldr$, the establishment of the covering condition for some open set $\cU\subset G$, is a criterion that implies the existence of orbits whose associated product of matrices in the fibers remain in $\cU$.

 \begin{theorem}\label{thm:covering for lc}
Let  $(\sigma,A)$ be a locally constant cocycle over the shift map $\sigma:\Sigma_m\to\Sigma_m$, depending on zeroth position and represented by $\phi:\cA\to G $.  Suppose that there exists a non-empty open set $\cU \subset G$ with compact closure satisfying
     \begin{equation}\label{eq:covering condition2}
          \overline{\cU}\subset \bigcup\limits_{i\in \cA} \phi(i)^{-1}\cU.
     \end{equation}
     Then, there exists 
 $\epsilon>0$, such that every $\alpha$-Hölder cocycle $(\sigma,B)$ with $d_{\mathcal{C}^{\alpha}}(A,B)<\epsilon$ such that 
 {there exists $x\in \Sigma_m$, such that for every $n\in \N$, $B_n(x)\in \cU$.}
 \end{theorem}
{
Theorem \ref{thm:covering for lc}, in particular, implies that  the covering condition \eqref{eq:covering condition2} for $(\sigma,A)$ implies that for every $B$ sufficiently $\cC^\alpha$-close to $A$, $(\sigma,B)\in \cB$.} A  generalization of the above theorem to general locally constant cocycles will be provided in 
Theorem \ref{thm:SFT-covering}.  
Before going into the proof of this theorem, we provide two lemmas regarding the control of distances between finite products of matrices associated with points in a cocycle that we need in the proof of Theorems \ref{thm:covering for lc} and  \ref{thm:SFT-covering}. 
    {\begin{lemma}\label{lem:finite-close}
     Let  $k\in \N$ and  $A,B:\Sigma_m\to \gldr$ be continuous maps. Then, for every $x\in \Sigma_m$,
     \begin{equation*}
      \|A_k(x)-B_k(x)\|  \leq k \max\{\|A\|,\|B\|\}^{k-1}d_{\cC^0}(A,B).   
     \end{equation*}
     \end{lemma}
     \begin{proof}
     For simplicity of notation, denote $A(\sigma^i(x))$ and $B(\sigma^i(x))$ by $A_i$ and $B_i$, respectively. Then, 
     \[A_k(x)-B_k(x)=\sum\limits_{i=1}^k A_k\cdots A_{i+1}(A_i-B_i)B_{i-1}\cdots B_1.\]
Employing the triangle inequality and using  $\|A_i-B_i\|\leq d_{\cC^0}(A,B)$, we get
\begin{align*}
  \|A_k(x)-B_k(x)\|& \leq \sum\limits_{i=1}^k \|A_k\|\cdots \|A_{i+1}\|.\|A_i-B_i\|.\|B_{i-1}\|\cdots \|B_1\| \\
  & \leq d_{\cC^0}(A,B) \sum\limits_{i=1}^k \|A\|^{k-i}\|B\|^{i-1}
\end{align*}
{
Now, for $C:=\max\{\|A\|, \|B\|\}$, we have $\sum\limits_{i=1}^k \|A\|^{k-i}\|B\|^{i-1}\leq k C^{k-1}$, completing the proof.} 
     \end{proof}

{
\begin{lemma}\label{lem:bound-Hölder}
Let $(\sigma,B)$ be an $\alpha$-Hölder $\gldr$ cocycle over the shift map. Assume that  there exist $L,m\in \N$, $K>0$ and a sequence $0=n_0<n_1<\cdots<n_k=m$ of  integers such that for every $x\in \Sigma$,  every $y\in \sigma^{-m}(W^u_{loc} \sigma^{m}(x))$ and every $i=1,\ldots,k$, 
\[n_{i}-n_{i-1}\leq L, \quad \|B_{n_i}(x)^{\pm 1}\|\leq K, \quad  \|B_{n_i}(y)^{\pm 1}\|\leq K.\]
Then,
\begin{equation}\label{eq:Hölder-better-bound}
    \|B_{m}(x)-B_{m}(y)\|<({1-2^{-\alpha}})^{-1}{K^{3}}\max\{\|B\|,1\}^{L-1}|B|_{\cC^\alpha}.
\end{equation}

\end{lemma}
\begin{proof}
To prove the Lemma, we first prove a prior upper bound for $\|B_m(x)-B_m(y)\|$ that works for every $x,y$, and for every $\alpha$-Hölder cocycle. Then, we use it as a part of the proof of the main statement. 

\medskip

\noindent
\textbf{Claim.} \textit{
Let $(\sigma,B)$ be an $\alpha$-Hölder cocycle over the shift map. 
Then,  for every $x,y$ 
\begin{equation}\label{eq:Hölder-ap-bound}
    \|B_m(x)-B_m(y)\|\leq \|B\|^{m-1}|B|_{\cC^{\alpha}}\sum\limits_{i=1}^m d(\sigma^i(x),\sigma^i(y))^\alpha.
\end{equation}
}

Denoting $x_i:=\sigma^{i}(x)$ and $y_i=\sigma^i(y)$, we have 
\[B_m(x)-B_m(y)=\sum\limits_{i=0}^{m-1} B_{m-i-1}(x_i)\big(B(x_i)-B(y_i)\big) B_i(x).\]
Therefore, by the triangle inequality 
\begin{align}
    \|B_m(x)-B_m(y)\| & \leq \sum\limits_{i=0}^{m-1} \|B_{m-i-1}(x_i)\|.\|B(x_i)-B(y_i)\|.\| B_i(x)\|\\
    & \leq \|B\|^{m-1}\sum\limits_{i=0}^{m-1} \|B(x_i)-B(y_i)\|\label{eq:first-bound}
\end{align}
By Hölder regularity of $B$, 
\(\|B(x_i)-B(y_i)\| \leq |B|_{\cC^{\alpha}}d(x_i,y_i)^\alpha\). Combining this with \eqref{eq:first-bound} gives \eqref{eq:Hölder-ap-bound} and finishes the proof of the claim. \\

Now, we turn to the proof of the main statement.  We can write 
      \begin{align*}\label{eq:expand}
     B_{n_k}(x)-B_{n_k}(y)&=\sum\limits_{i=0}^{k-1} B_{n_k-n_{i+1}}(x_{n_{i+1}})\Big(B_{n_{i+1}-n_{i}}(x_{n_{i}})-B_{n_{i+1}-n_{i}}(y_{n_{i}})\Big)B_{n_{i}}(y).
\end{align*}
Note that since 
\(B_{n_k-n_{i+1}}(x_{n_{i+1}})=B_{n_k}(x)B_{n_{i+1}}(x)^{-1}\), by assumption, $\|B_{n_i}(x)^{\pm 1}\|\leq K$, we get $\|B_{n_k-n_{i+1}}(x_{n_{i+1}})\|\leq K^2$. Therefore, applying the triangle inequality to the above equality, we get 
    \begin{align}
        \| B_{n_k}(x)-B_{n_k}(y)\| \leq K^{3}\sum\limits_{i=0}^{k-1} \|B_{n_{i+1}-n_{i}}(x_{n_i})-B_{n_{i+1}-n_{i}}(y_{n_i})\|.
    \end{align}
    On the other hand, by \eqref{eq:Hölder-ap-bound} and $n_{i+1}-n_i\leq L$, we have for every $i$
    \begin{align}
        \|B_{n_{i+1}-n_{i}}(x_{n_i})-B_{n_{i+1}-n_{i}}(y_{n_i})\| & \leq \|B\|^{n_{i+1}-n_i-1}|B|_{\cC^\alpha} \sum\limits_{j=n_{i-1}+1}^{n_i} d(x_j,y_j)^\alpha \\   
        & \leq \max\{\|B\|,1\}^{L-1}|B|_{\cC^\alpha} \sum\limits_{j=n_{i-1}+1}^{n_i} d(x_j,y_j)^\alpha.
    \end{align}

Therefore, 
\begin{equation}
     \| B_{n_k}(x)-B_{n_k}(y)\| \leq K^3 \max\{\|B\|,1\}^{L-1}|B|_{\cC^\alpha}\sum\limits_{i=0}^{m-1}d(x_i,y_i)^\alpha.
\end{equation}
Now, since $y\in \sigma^{-m}(W^u_{loc}\sigma^m(x))$, for every  $i=0,\ldots,m-1$, $d(x_j,y_j)^\alpha
     < (\frac{1}{2^{m-j}})^\alpha|B|_{\cC^\alpha}$, which implies
\begin{align*}
    \|B_{m}(x)-B_{m}(y)\| &  \leq  K^{3}\|B\|^{L-1}|B|_{\cC^\alpha} \sum\limits_{i=0}^{m-1} (\frac{1}{2^{m-i}})^\alpha \\ & \leq  K^{3}\max\{\|B\|,1\}^{L-1}|B|_{\cC^\alpha}\Big(\sum\limits_{i=0}^{\infty}\frac{1}{2^{i\alpha}}\Big) \\
    & = ({1-2^{-\alpha}})^{-1}{K^{3}}\max\{\|B\|,1\}^{L-1}|B|_{\cC^\alpha}.
    \end{align*}
This completes the proof. 
\end{proof} 
}

     {\begin{proposition}\label{prop:finite-step-qc}
         Under assumptions of Theorem \ref{thm:covering for lc}, there exists  $\epsilon>0$ such that
         for every $\alpha$-Hölder cocycle $(\sigma,B)$ with $d_{\mathcal{C}^{\alpha}}(A,B)<\epsilon$, there is a sequence $\{x_n\}_{n\geq 0}$ in $\Sigma_m$ such that for every $n\geq 0$, if  $W_n:= W^u_{loc}(\sigma^n(x_n))$, then 
         \begin{itemize}
             \item[(i)] $\sigma^{-(n+1)}(W_{n+1})\subset \sigma^{-n}(W_n)$,
             \item[(ii)]  for every $i=1,\ldots,n+1$ and $y\in \sigma^{-n}(W_n)$,  $B_i(y) \in \cU$.
         \end{itemize}
     \end{proposition}}
     
     \begin{proof}
     
        The proof is by induction.
        {For every $j\in \cA$, denote $D_j:=\phi(j)$. First, in view of {Lemma \ref{lem:smaller-U}},  we may assume that $\cU$ contains the identity matrix, $\id$. Additionally, there is $\delta>0$ such that 
        \begin{align} \label{eq:covering with delta}
             \overline{\cU}\subset \bigcup\limits_{i=1}^m {D_i}^{-1}\cU_{(\delta)}.
         \end{align}
         }

    Since $\id\in \cU$, it follows from \eqref{eq:covering with delta} that there is $i_0\in \cA$ such that $D_{i_0}\in \cU_{(\delta)}$. Clearly, $A|_{\cylinder[i_0]}=D_{i_0}$. Choose an element of ${\cylinder[i_0]}$ and denote it by $x_0$. Since $W_0=W^u_{loc}(x_0) \subset {\cylinder[i_0]}$, we conclude that  for every $y\in W_0$, $A(y)=D_{i_0}\in \cU_{(\delta)}$.

Now, if $\epsilon\in (0,\delta)$ and  $d_{\mathcal{C}^{\alpha}}(A,B)<\epsilon$, then $d_{\mathcal{C}^{0}}(A,B)<\epsilon$. Consequently,  for every $y\in W_0$, 
$B(y)\in \cU$.

The above argument proves the proposition for $n=0$. Now, suppose that the proposition holds for $n-1$, that is there exists $x_{n-1}\in \Sigma_m$ satisfying both items of the proposition.
We aim to find $x_{n}\in \Sigma_m$ satisfying the desired properties. Note that by \eqref{eq:covering with delta}, 
\[B_{n}(x_{n-1})\in \overline{\cU}\subset \bigcup\limits_{i=1}^m {D_i}^{-1}\cU_{(\delta)}.\]
Therefore, there exists $i_{n}\in\cA$ such that  $D_{i_{n}}B_{n}(x_{n-1})\in \cU_{(\delta)}$ and since $A$ is locally constant depending on the zeroth position, $A|_{\cylinder[i_n]}=D_{i_{n}}$. Let $x_{n}$ be an arbitrary element in $\sigma^{-(n-1)}(W_{n-1})\cap \sigma^{-n}(\cylinder[i_n])$. We claim that this choice of $x_{n}$ satisfies the required properties. (See Figure \ref{fig:nested-u-segments}).

\medskip

\begin{figure}
    \centering
\includegraphics[width=\textwidth]{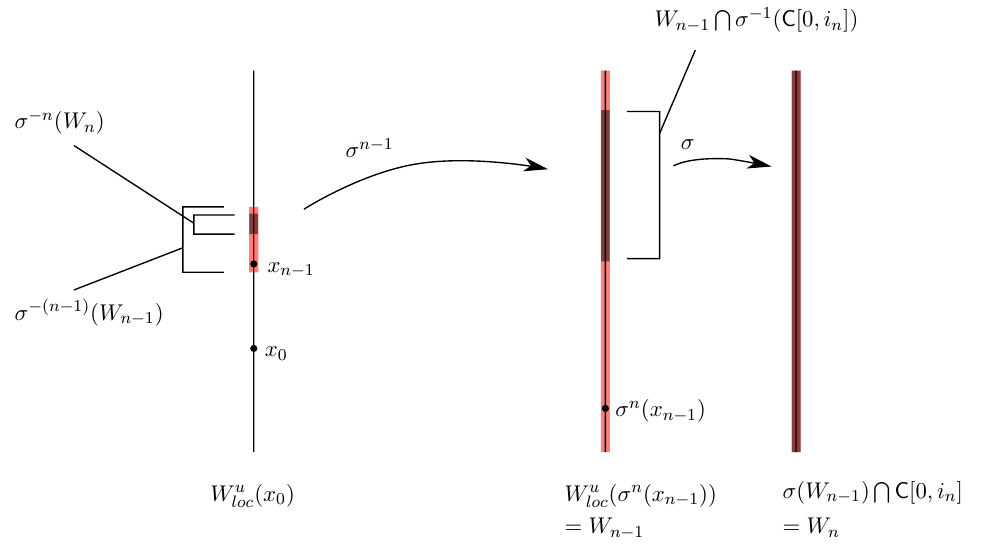}
    \caption{Proof of Proposition \ref{prop:finite-step-qc}. }
    \label{fig:nested-u-segments}
\end{figure}

\noindent
\textbf{Claim}. \textit{There exists $C>0$ depending on $\cU,A,\alpha$ such that for every $z\in \sigma^{-(n+1)}(W_n)$,} \[\|B_{n+1}(z)-D_{i_n}B_{n}(y)\|<C\epsilon.\]

\medskip
To prove the claim, first using the triangle inequality, one has
\begin{equation}\label{equ1}
  \|B_{n+1}(z)-D_{i_n}B_{n}(y)\|\leq  \|B_{n+1}(z)-D_{i_n}B_{n}(z)\|+ \|D_{i_n}B_{n}(z)-D_{i_n}B_{n}(y)\|. 
\end{equation}
Since $D_{i_n}=A|_{\cylinder[i_n]}$ and $\sigma^{n}(z)\in \cylinder[i_n]$, and $d_{\cC^{\alpha}}(A,B)<\epsilon$, $\|B(\sigma^{n}(z))-D_{i_n}\|<\epsilon$. 
On the other hand, $z\in \sigma^{-(n+1)}(W_{n+1}) \subset \sigma^{-n}(W_n)$. Therefore,
\begin{equation}\label{equ3}
    \|B_{n+1}(z)-D_{i_n}B_{n}(z)\| \leq \|B(\sigma^{n}(z))-D_{i_n}\|.\|B_{n}(z)\|\leq  K_{\cU}\epsilon,
\end{equation}
where $K_{\cU}:=\max\{\|D^{\pm 1 }\|:D\in \overline{\cU}\}$. 

To estimate the second term in \eqref{equ1}, using Lemma \ref{lem:bound-Hölder}, we get
\begin{align}
     \|D_{i_n}B_{n}(y)-D_{i_n}B_{n}(z)\| & \leq \|D_{i_n}\| . \|B_{n}(y)-B_{n}(z)\|\\
     & \leq (1-2^{-\alpha})^{-1}\|A\|K_{\cU}^{3}|B|_{\cC^\alpha}.\label{difference-upper-bound}
     \end{align}

     Since all the points $y,z$ and $x_i$'s are in the same local unstable set, we can replace $|B|_{\cC^\alpha}$ by $\sup_{x \in \Sigma}|B\big\vert_{W^u_{loc}(x)}|_{\cC^\alpha}$ in (\ref{difference-upper-bound}). On the other hand, since  $A$ is a locally constant cocycle depending on zeroth position,
it follows from   Corollary \ref{cor:Hölder}, $\sup_{x \in \Sigma}|B\big\vert_{W^u_{loc}(x)}|_{\cC^\alpha}<\epsilon$. So 
     \begin{equation}\label{equ2}
          \|D_{i_n}B_{n}(y)-D_{i_n}B_{n}(z)\|  \leq (1-2^{-\alpha})^{-1}\|A\|K_{\cU}^{3}\epsilon.
     \end{equation}
 The estimates provided by (\ref{equ1}), (\ref{equ3}) and (\ref{equ2}), complete the proof of Claim for $C=K_{\cU}+(1-2^{-\alpha})^{-1}\|A\|K_{\cU}^{3}$. 

To complete the proof of the proposition, it is enough to assume that 
$\epsilon<C^{-1}\delta$, for $C$ given by the claim.  Then, according to the claim, for every $z\in \sigma^{-(n+1)}(W_n)$,  $ \|B_{n+1}(z)-D_{i_n}B_{n}(y)\| \leq \delta$ and since $D_{i_n}B_{n}(y)\in \cU_{(\delta)}$, we get $B_{n+1}(z)\in \cU$.
     \end{proof}

Now we can prove Theorem \ref{thm:covering for lc}.

\begin{proof}[Proof of Theorem \ref{thm:covering for lc}]
{
Let $\epsilon>0$ be the number provided by Proposition \ref{prop:finite-step-qc}. Take cocycle $B$ with $d_{\mathcal{C}^{\alpha}}(A,B)<\epsilon$ and  $\{x_n\}_{n\in \N}$ be the sequence given in Proposition \ref{prop:finite-step-qc}.  
According to this proposition, for $W_n:=W^u_{loc}(\sigma^n(x_n))$, $\{\sigma^{-n}(W_n)\}_{n\in \N}$ forms a nested sequence of compact subsets of $\Sigma_m$, ensuring that their intersection is non-empty. Now, in view of  Proposition \ref{prop:finite-step-qc}, for every point $x^*$ in this intersection $B_n(x^*)\in \cU$ for every $n\geq 0$.
}
\end{proof}
 
 \subsection{General locally constant cocycles}\label{subsec:stability-qc}

In this subsection, we concentrate on general locally constant cocycles and show that the covering criterion \eqref{eq:gen-covering} ensures the presence of quasi-conformal orbits and their stability in the $\cC^\alpha$ topology.

 \begin{theorem}[Key theorem]\label{thm:SFT-covering}

 {Let $G$ be a {closed} subgroup of $\gldr$.}
     Let $(\sigma,A)$ be a locally constant $G$ cocycle over the full shift on finitely many symbols,  depending on positions in $\ldbrack -k,k\rdbrack$ and represented by $\phi:\cA^{2k+1}\to G$.
     
      {
     If $(\sigma,A)$ satisfies the covering condition (for some non-empty open set $\cU\subset G$ with compact closure), then there exists $\epsilon,K>0$, such that  for every  $\alpha$-Hölder map $B:\Sigma_m\to G$ with $d_{\mathcal{C}^{\alpha}}(A,B)<\epsilon$,  the cocycle $(\sigma,B)$ admits $x\in \Sigma_m$ with \[\sup\limits_{n\in \Z}\|B_n(x)^{\pm 1}\|\leq K.\]
     In particular, $(\sigma,B)\in \cB$. 
     } 
 \end{theorem}

 The following proposition provides the main ingredient of the proof of Theorem \ref{thm:SFT-covering}. This is an analog of Proposition \ref{prop:finite-step-qc} for locally constant cocycles depending on positions in $\ldbrack -k,k\rdbrack$ for $k\geq 1$. 
 
\begin{proposition}
\label{prop:stability-Hölder-calculation}

  {Let $G$, $(\sigma,A)$ and $\cU$ satisfy the assumptions of Theorem \ref{thm:SFT-covering}.}
  Then, there is  $\epsilon>0$ such that for every $\alpha$-Hölder {map $B:\Sigma_m\to G$} with $d_{\mathcal{C}^\alpha}(A,B)<\epsilon$, there exist 
     \begin{itemize}
        \item a positive integer $L$,
         \item a sequence ${0=n_0<}n_1<n_2<\cdots$ of integers,
         \item a sequence $\{x_{j}\}_{j\in \N}$ in $\Sigma_m$, 
         
     \end{itemize}
       such that for $W_j:=W^u_{loc}\big(\sigma^{n_{j}}({x_{j}})\big)$ and every $j\geq 0$, 
\begin{itemize}
    \item[(i)]  $n_{j+1}-n_j\leq L$,
    \item[(ii)] $\sigma^{-n_{j+1}}(W_{j+1}) \subset \sigma^{-n_j}(W_j)$, and
    \item[(iii)]     for every $y\in \sigma^{-n_j}(W_j)$ and $i=0,\ldots,j$, $B_{n_i}(y)\in\cU$.
\end{itemize}
\end{proposition}

\begin{proof}{
Suppose that the cocycle $(\sigma,A)$ satisfies the covering condition for  
$\cU\subset G$  and
$\cA'\subset \cA^{\ldbrack -k,k\rdbrack}$. In view of Remark \ref{rmk:better-gen-cov}, we may assume that $\cU$ contains the identity matrix and moreover, there exist $\delta>0$ and finite families $\Lambda_{\ua}$ ($\ua\in \cA'$) of transitions, such that for every $\ua\in \cA'$,  $\Lambda_{\ua}\subset \bigcup_{\ub\in \cA'}T(\ua,\ub)$ and
 \begin{equation}\label{eq:smaller-u-t}
      \overline{\cU} \subset \bigcup\limits_{\uw\in \Lambda_{\ua}}A_{[\uw\rangle}^{-1}\cU_{(\delta)}.
\end{equation}

Now, we take $L:=\max\{|\uw| :  \uw\in \Lambda_{\ua},\ua\in \cA'\}$ and  construct the sequence $\{n_j\}_{j\geq 0}$ so  that $n_0=0$ and for every $j\geq 0$, $0<n_{j+1}-n_j\leq L$. 
 First, fix an arbitrary element $\ua_0 \in \cA'$,  and set $x_0$ as an arbitrary element in $\cylinder[-k;\ua_0]$.  

\medskip

\noindent
 \emph{{\bf Claim.} There exist $\epsilon>0$ and   $\ua_1,\ua_2,\ldots$ in $\cA'$ such that for every $\alpha$-Hölder cocycle $(\sigma,B)$ with $d_{\cC^\alpha}(A,B)<\epsilon$, every $j\in \N$,   
 \[B_{n_j}(y)\in \cU, \quad \text{for every~}y\in \cylinder[-k;\uw_1\vee\cdots \vee \uw_j]\cap W^u_{loc}(x_0),\] 
 where $ \uw_j \in T(\ua_{j-1},\ua_j)$ and $n_j:=-k+\sum_{i\leq j}|\uw_i|$.}

 \medskip
 
    We prove the claim inductively. Assume that $\ua_0,\ldots,\ua_j\in \cA'$  have been chosen satisfying the conditions.  Let $x_j$ be an arbitrary point in $\cylinder[-k;\uw_1\vee \cdots \vee \uw_j]\cap W^u_{loc}(x_0)$. Then,  by \eqref{eq:smaller-u-t}, 
    \[B_{n_j}(x_j)\in \overline{\cU} \subset \bigcup\limits_{\uw\in \Lambda_{\ua_j}}A_{[\uw\rangle}^{-1}\cU_{(\delta)}.\]
    Therefore, there exists $\ua_{j+1}\in \cA'$  and $\uw_{j+1}\in T(\ua_j,\ua_{j+1})$ such that 
    \begin{align}\label{covering-induction}
     A_{[\uw_{j+1}\rangle}B_{n_j}(x_j)\in \cU_{(\delta)}.
    \end{align}
    Next, we aim to show that for every point $z\in \cylinder[-k; \uw_1\vee \cdots \vee \uw_{j+1}]\cap W^u_{loc}(x_0)$, $B_{n_{j+1}}(z)\in \cU$, where $n_{j+1}=n_j+|\uw_{j+1}|$. To see this, first observe that 
\begin{align*}
    B_{n_{j+1}}(z)&-A_{[\uw_{j+1}\rangle}B_{n_j}(x_j)\\
    &  = \big(B_{n_{j+1}-n_j+k}(\sigma^{n_j-k}(z))-A_{[\uw_{j+1}\rangle}B_k(\sigma^{n_j-k}(x_j))\big)B_{n_j-k}(z)\\ 
  &  \; +A_{[\uw_{j+1}\rangle}B_k(\sigma^{n_j-k}(x_j))\big(B_{n_j-k}(z)-B_{n_j-k}(x_j)\big).
\end{align*}
Therefore, by the triangle inequality, 
\begin{align}\label{eq:triangle-ineq}
\begin{split}
  \|  B_{n_{j+1}}(z)-& A_{[\uw_{j+1}\rangle}B_{n_j}(x_j)\|  \\
  & \leq \|B_{n_{j+1}-n_j+k}(\sigma^{n_j-k}(z))-A_{[\uw_{j+1}\rangle}B_k(\sigma^{n_j-k}(x_j))\|.\|B_{n_j-k}(z)\| \\
  & + \|A_{[\uw_{j+1}\rangle}B_k(\sigma^{n_j-k}(x_j))\|. \|B_{n_j-k}(z)-B_{n_j-k}(x_j)\|.
\end{split}
\end{align}
Note that since $\cylinder[-k; \uw_1\vee \cdots \vee \uw_{j+1}]\subset \cylinder[-k; \uw_1\vee \cdots \vee \uw_{j}]$, it follows from induction hypothesis that $B_{n_j}(z)\in \cU$ and so 
\begin{equation*}
    \|B_{n_j}(z)\|\leq K_{\cU}:=\sup\{\|D^{\pm 1}\| :  D\in \overline{\cU}\}. 
\end{equation*}
Therefore, 
\begin{equation*}
    \|B_{n_j-k}(z)\|\leq \|B_{n_j}(z)\|.\|B_k(\sigma^{n_j-k}(z))^{-1}\|\leq K_{\cU}\|B^{-1}\|^k\leq K_{\cU}\|B\|^{kd}.
\end{equation*}
Now, since $d_{\cC^0}(A,B)< \epsilon$, one has  $\|B\|\leq \|A\|+\epsilon$, and so 
\begin{equation}\label{eq:bound-Bnj-k'}
    \|B_{n_j-k}(z)\|\leq K_{\cU}(\|A\|+\epsilon)^{kd}.
\end{equation}
On the other hand,  
it follows from Lemma
 \ref{lem:finite-close} that 
\begin{align}\label{eq:boundBnj+1}
 \begin{split}
  \|B_{n_{j+1}-n_j+k}(\sigma^{n_j-k}(z))-A_{[\uw_{j+1}\rangle}&B_k(\sigma^{n_j-k}(x_j))\|\\ &\leq (L+k)(\|A\|+\epsilon)^{L+k-1}\epsilon.
 \end{split}
\end{align}

Meanwhile, by Lemma \ref{lem:bound-Hölder}, we have
\begin{align}
       \|  B_{n_{j}-k}(z)-B_{n_j-k}(x_j)\| &\leq (1-2^{-\alpha})^{-1}\|B\|^{L-1}{K_{\cU}^{3}}|B|_{\cC^\alpha}, \label{eq:bound-1}
\end{align}

Since all the points $z$ and $x_j$'s and their relevant iterations are in the same local unstable set with distance less than $2^{-k}$, we can replace $|B|_{\cC^\alpha}$ in \eqref{eq:bound-1} by $\sup_{x \in \Sigma}|B\big\vert_{\sigma^{-k}(W^u_{loc}(x))}|_{\cC^\alpha}$. On the other hand, since  $A$ is a locally constant cocycle depending on positions in $\ldbrack -k,k\rdbrack$,
it follows from   Corollary \ref{cor:Hölder}, 
\begin{equation}
    \sup_{x \in \Sigma}|B\big\vert_{\sigma^{-k}(W^u_{loc}(x))}|_{\cC^\alpha}<\epsilon.
\end{equation}
Therefore, by \eqref{eq:bound-1},
\begin{align}
       \|  B_{n_{j}-k}(z)-B_{n_j-k}(x_j)\| &\leq (1-2^{-\alpha})^{-1}\|B\|^{L-1}{K_{\cU}^{3}}\epsilon.\label{eq:bound-1'}
\end{align}

Now,
by 
\eqref{eq:triangle-ineq}--\eqref{eq:bound-1'},
we have 
\begin{equation}\label{eq:final-bound}
     \|  B_{n_{j+1}}(z)-A_{[\uw_{j+1}\rangle}B_{n_j}(x_j)\|\leq C\epsilon.
\end{equation}
for some constant $C$ depending on $\cU,L,d,\alpha,\|A\|,k$ and independent of $j$. So, due to  \eqref{covering-induction} and \eqref{eq:final-bound}, if $\epsilon<C^{-1}\delta$, we deduce that $B_{n_{j+1}}(z)\in \cU$ and the proof of Claim is complete. }\end{proof}

\begin{proof}[Proof of Theorem \ref{thm:SFT-covering}] 
{
Let $(\sigma,A)$ be a cocycle that satisfies the given assumptions. According to Proposition \ref{prop:stability-Hölder-calculation}, there exists $\epsilon > 0$ such that for every $\alpha$-Hölder cocycle $(B, \sigma)$ with $d_{\mathcal{C}^{\alpha}}(A, B) < \epsilon$, we can find a sequence $0 = n_0 < n_1 < n_2 < \cdots$ of integers with $L := \sup_j (n_{j+1} - n_j) < +\infty$, a sequence $\{x_j\}_{j \geq 1}$ in $\Sigma_m$, and a non-empty open set $\mathcal{U} \subset G$ with compact closure such that for $W_j = W^u_{loc} \big(\sigma^{n_j}(x_j)\big)$, $\sigma^{-n_{j+1}}(W_{j+1}) \subset \sigma^{-n_{j}}(W_j)$,  and for every $y \in W_j$, $B_{n_j}(y) \in \mathcal{U}$.

Since $\{\sigma^{-n_{j}}(W_j)\}_{j \geq 0}$ forms a nested sequence of compact subsets of $\Sigma_m$, their intersection must be non-empty. We claim that for every point $x^*$ in this intersection, the sequence $\{\|B_n(x^*)\|\}_{n \in \mathbb{N}}$ is bounded by some constant depending on $A$.

If $x^*$ belongs to this intersection, then for every $j \in \mathbb{N}$, $B_{n_j}(x^*) \in \mathcal{U}$ and so 
\[ \|B_{n_j}(x^*)^{\pm 1}\| \leq K_\cU,\] 
where $K_\cU := \sup\{\|D^{\pm 1}\| : D \in \overline{\mathcal{U}}\}$.
Given $n \in \mathbb{N}$, there exists $j \geq 0$ such that $n_j \leq n < n_{j+1} \leq n_j + L$. We use the identity $B_n(x^*)=B_{n-n_j}(\sigma^{n_j}(x^*))B_{n_j}(x^*)$ to give upper bound for $\|B_n(x^*)^{\pm 1}\|$. We have,
\begin{align*}
    \|B_n(x^*)\| & \leq \|B_{n-n_j}(\sigma^{n_j}(x^*))\|.\|B_{n_j}(x^*)\| \leq \|B\|^{n-n_{j}} K_\mathcal{U}.
\end{align*} 
Similarly, $ \|B_n(x^*)^{-1}\|\leq \|B^{-1}\|^{n-n_{j}} K_\mathcal{U}$.
 Now, since  $\|A-B\|<\epsilon$,  we get that if $\epsilon$ is sufficiently small, there exists $C$ depending on $\epsilon$ and on cocycle $A$ such that $\|B^{-1}\|< \|A^{-1}\|+C\epsilon$. Hence,
 \[\|B_n(x^*)^{\pm 1}\|\leq  C'K_\cU,\]
 for some constant $C'$ depending on $\epsilon,L$ and the cocycle $A$. 
 Finally, using  Lemma \ref{lem:bounded-properties} part (iv), we conclude that for every $y^*\in \omega(x^*)$, 
  \[\sup_{n\in \Z}\|B_n(y^*)^{\pm 1}\| \leq  K:=(C'K_\cU)^2.\]
 This completes the proof.
} \end{proof}

 \subsection{The case of general linear group}\label{subsec:GL-qc}

{As explained at the beginning of the section, Theorem \ref{thm:covering for lc} holds for every subgroup of $\gldr$. Let $G=\gldr$. Hence, to prove the Theorem \ref{thm:GL-bd-example}, we must find a finite family of matrices in $\gldr$ to satisfy the covering condition for some open subset of $\gldr$.}
 \begin{proof}[Proof of Theorem \ref{thm:GL-bd-example}] 
{Using Corollary \ref{cor:gldr-covering}, there exist $\{D_1,\ldots,D_{2d^2}\}$  of matrices in $\gldr$ and open set $\cU \subset \gldr$ which satisfy the covering condition. Let $\sigma:\Sigma_{2d^2} \to \Sigma_{2d^2}$ be the shift map. Define a locally constant cocycle $(\sigma,A)$, depending on zeroth position and associated to the IFS generated by $\langle D_1,\ldots,D_{2d^2}\rangle^+$. Thus Theorem \ref{thm:covering for lc} proves that, there exists $\epsilon>0$ such that  if $d_{\mathcal{C}^{\alpha}}(A,B)<\epsilon$, for a $\gldr$ cocycle $(\sigma,B)$, then there exists $x^* \in \Sigma_{2d^2}$, such that for every $n\in \N$
\[\|B_n(x^*)^{\pm 1}\|\leq \kappa_{\cU}.\]
 To conclude the proof, notice that by Lemma \ref{lem:bounded-properties} part (iv), for every point $y^*\in \omega(x^*)$, $\omega$-limit set of orbit of $x^*$, we have  
 \(\|B_n(y^*)^{\pm 1}\| \leq \kappa^2_{\cU}\) and this completes the proof.}
\end{proof}

\section{Cocycles without domination after {\cite{BDP}}}\label{sec:transition}

In this section, we consider the class of cocycles called cocycles with transitions.  It was introduced in \cite{BDP} and played
a crucial for the perturbative arguments needed in the proofs.

\subsection{Transition property}\label{subsec:transition}
A  \textit{periodic cocycle} is a cocycle such that the base dynamics contain only periodic orbits. We denote by $n(x)$ the period of the periodic point $x$.

\begin{definition} {\cite[Definition 1.6]{BDP}}
   Given $\epsilon>0$, a periodic cocycle $(\boldsymbol{E},X,F,f)$ admits \textit{$\epsilon$-transitions}, if for every finite family of points $x_1,\ldots,x_n \in X$, there is a local trivialization or an orthonormal system of coordinates of the vector bundle $(\boldsymbol{E},X,\pi)$ (now fix some matrix representation), and for any $(i,j) \in \{1,\ldots,n\}^2$ there exists a matrix transition $[T_{i,j}]$ from $[F_i]:=\big(F(f^{n(x_i)-1}(x_i)),\ldots,F(x_i)\big)$ to $[F_j]:=(F(f^{n(x_j)-1}(x_j)),\ldots,F(x_j))$   satisfying the following properties:
   
   \begin{itemize}
       \item[1)]  For every $m \in \N$, $\underline{\ell}=({\ell}_1,\ldots,{\ell}_m) \in \{1,\ldots,n\}^m$, and $\ua=(a_1,\ldots,a_m) \in \N^m$ consider the word
   $$[W(\underline{\ell},\ua)]=[T_{\ell_1,\ell_m}][F_m]^{a_m}[T_{\ell_m,\ell_{m-1}}] \cdots [T_{\ell_2,\ell_1}][F_1]^{a_1},$$
   where  the word $((x_{\ell_1},a_1),\ldots,(x_{\ell_m},a_m))$ with letters in $X\times \N$ is not a power of a shorter word. 
  
 Then, there is $x\in X$ of period equal to the length of $[W({\ell,\ua})]$ such that:
   \begin{itemize}
      \item   the word $(F(f^{n(x)-1}(x)),\ldots,F(x))$ is  $\epsilon$-close\footnote{Two words $[T]=(T_1,\ldots,T_n)$ and $[T']=(T'_1,\ldots,T'_n)$  of matrices with the same length are called $\epsilon$-close if $\max\limits_{1\leq i \leq n}\|T_i-T'_i\|<\epsilon$.} to $[W(\underline{\ell},\ua)]$,
      \item there exists an $\epsilon$-perturbation $\Tilde{F} $ of $F$ in the $\cC^0$ topology  such that $$[W(\underline{\ell},\ua)]=(\tilde{F}(f^{n(x)-1}(x)),\ldots,\tilde{F}(x)).$$      
   \end{itemize}

\item[2)] One can choose $x$ such that the distance between the orbit of $x$ and any point $x_{\ell_i}$ be arbitrary small provided that $a_i$ is sufficiently large.
\end{itemize}
   \end{definition}

The following is compatible with the one in \cite{BDP}.
   \begin{definition}\label{def:gen-transition}
   A cocycle $(\boldsymbol{E},X,F,f)$ has {\it transitions property} if the periodic points of $f$ are dense in $X$ and the restriction of $f$ to its periodic points admits $\epsilon$-transition for every $\epsilon>0$.
\end{definition}

\begin{example}
    By Lemma \ref{lem: locall-realization}, every locally constant cocycle has transitions property. Indeed, 
Lemma \ref{lem: locall-realization} allows us to 
realize a product of matrices over a periodic fiber for a locally constant cocycle using a transition. Generally speaking, a cocycle has transitions property if it can realize the product of matrices after arbitrary small perturbations.
\end{example}

    \begin{lemma}
    Let $\sigma:\Sigma_H\to \Sigma_H$ be a transitive subshift of finite type, then every $\sldr$ cocycle $(\sigma,A)$ over $\sigma:\Sigma_H\to \Sigma_H$  has transitions property. 
\end{lemma}

\begin{proof} 
 It is known that transitivity of $\sigma:\Sigma_H\to \Sigma_H$ implies that the set of periodic points is dense in $\Sigma_H$. So, we need to verify the transition property. 
 
Consider two periodic points $p,p'\in \Sigma_H$ of periods $n,n'$, respectively.  
Let $\ua,\ua'\in \cA^{\mathrm{fin}}$ with $|\ua|=n$ and $|\ua'|=n'$ be such that $p=(\ldots,\ua;\ua,\ua,\ldots)$ and  $p'=(\ldots,\ua';\ua',\ua',\ldots)$. 
Since $\sigma:\Sigma_H\to \Sigma_H$ is transitive, there exists $\ub,\underline{c}\in \cA^{\mathrm{fin}}$  such that the concatenated words $\ua'\uc,\uc \,\ua,\ua \, \ub,\ub \, \ua'$ are all admissible (with respect to $H$). This implies that for every $i,j\in \N$, the  periodic point for $\sigma$ associated to the following finite word is an element of $\Sigma_H$:
\[\uw_{ij}:=\uc\underbrace{\ua\,\ua\cdots \ua}_{i}\ub\underbrace{\ua'\ua'\cdots \ua'}_{j}.\]

Now, fixing $\epsilon>0$,
since $A:\Sigma_H\to \sldr$ is uniformly continuous,  there exists  $k \in \N$ such that $\|A(x)-A(y)\|<\epsilon$ whenever $x,y\in \Sigma_H $ with $d(x,y)<2^{-k}$.

Take $i,j>\max\{2k/n,2k/n'\}$ and denote the periodic point associated to $\uw_{ij}$ by  $p^*$. Then, it is easy to check that the following is an $\epsilon$-transition from $p'$ to $p$
{\begin{equation*}
\Big(A(\sigma^{-kj}(p^*)),A(\sigma^{-kj+1}(p^*)),\ldots,A(p^*),A(\sigma(p^*)),\ldots,A(\sigma^{|\uc|+ki}(p^*))\Big),
\end{equation*}
as required.}
\end{proof}

\subsection{Cocycles without dominated splitting}\label{subsec:no-domination}

{
\begin{proposition}[{\cite[Proposition 2.1]{BDP}}]\label{prop:BDP-dense-dichotomy}
Let $d\geq 2$ and $(\boldsymbol{E},X,F,f)$ be a $d$-dimensional periodic cocycle with transitions which does not admit dominated splitting. Then, for every $\epsilon>0$,  there exist an $\epsilon$-perturbation $\tilde{F}$ of $F$  in the $\cC^0$ topology and a point $x\in X$ such that $\tilde{F}_{n(x)}(x)$ is a homothety. 
\end{proposition}}
{
A main part in the proof of Proposition \ref{prop:BDP-dense-dichotomy} in \cite{BDP} is a dichotomy in the family of linear cocycles with transitions (See \cite[Proposition 2.4]{BDP}). For every $\epsilon>0$, $d\geq 2$ and every $i=1,2,\ldots,d-1$,  every $d$-dimensional cocycle $(\boldsymbol{E},X,F,f)$ with transitions:
\begin{itemize}
    \item[(i)]either $(\boldsymbol{E},X,F,f)$  admits a dominated splitting $\boldsymbol{E}_2\prec \boldsymbol{E}_1$ {of index $i$.}
    \item[(ii)] or there exists an $\epsilon$-perturbation $\tilde{F}$ of $F$ and a periodic point $p\in X$ such that denoting the spectrum of  $\tilde{F}_{n(p)}(p)$ in increasing absolute value by $\lambda_1,\ldots,\lambda_d$, then $\lambda_i,\lambda_{i+1}$ are non-real complex conjugate numbers and  if $1<i$, then $\log |\lambda_{i-1}|<\log |\lambda_{i}|$ and if $i<d-1$, then $\log |\lambda_{i+1}|<\log |\lambda_{i+2}|$.  

\end{itemize}

 We call a periodic point such as $p$ in the second item above an $i$-elliptic periodic point (cf. \cite{Gourmelon2022}). One can easily verify that the existence of an $i$-elliptic periodic point is a $\cC^0$-stable property, so if a cocycle with transitions does not admit any dominated splitting, then there exists some perturbation of cocycle which has $i$-elliptic periodic points for every $i=1,2,\ldots,d-1$. Therefore, we can conclude the following corollary.
\begin{corollary}\label{cor:stably no dom}
    Let $\epsilon>0$ and $(\boldsymbol{E},X,F,f)$ be a cocycle with transitions that does not admit any dominated splitting. Then, there exists an $\epsilon$-perturbation  $\tilde{F}$ of $F$ in the  $\cC^0$ topology, such that $(\boldsymbol{E},X,\tilde{F},f)$ stably does not admit any dominated splitting.
\end{corollary}}

{Since every dominated splitting on a set of periodic points can be extended to a dominated splitting over the closure of periodic points \cite[B.1]{BDV} and referring to the Definition \ref{def:gen-transition}, we can consider cocycles with transitions instead of periodic cocycles with transitions.}
By normalizing the cocycles, one obtains the following variant of this result for the volume-preserving cocycles. 

\begin{corollary}\label{cor:BDP-sl}
   For any $d\geq 2$ and $\epsilon>0$ there is $l>0$ such that any $d$-dimensional volume preserving cocycle $(\boldsymbol{E},X,F,f)$ with transitions, satisfies the following:
\begin{itemize}
    \item either $(\boldsymbol{E},X,F,f)$  admits a {dominated splitting of index $l$},
    \item or there exist an $\epsilon$-perturbation $\tilde{F}$ of $F$ in the $\cC^0$ topology and a periodic point $x\in X$ such that $\tilde{F}_{n(x)}(x)=\id $.
    \end{itemize}
\end{corollary}

 We need the following refinement of the results above, which is the aim of this section.

 \begin{theorem} \label{thm:BDP-multi-orbits}
 Let $d\geq 2$. 
 Suppose that $(\boldsymbol{E},X,F,f)$ is a $d$-dimensional volume-preserving cocycle with transitions that stably does not admit any dominated splitting. {Then, for every $\epsilon>0$ and $N\in \N$ there is a volume-preserving cocycle $\Tilde{F}$ over $f$ and  periodic points $p_1,\ldots,p_N$ for $f$ with disjoint orbits such that $d_{\cC^0}(F,\tilde{F})<\epsilon$  and for every $1\leq i \leq N$, $\tilde{F}_{n(p_i)}(p_i)=\id $.}
 \end{theorem}

  \begin{proof}[Proof of Theorem \ref{thm:BDP-multi-orbits}]
The assumption that the cocycle $(\boldsymbol{E},X,F,f)$ stably does not admit dominated splitting guarantees that there exists $\delta_0>0$ such that for every continuous volume-preserving cocycle $\hat{F}$ over $f$, if  
 $d_{\cC^0}(F,\hat{F})<\delta_0$, then $(\boldsymbol{E},X,\hat{F},f)$ does not admit dominated splitting.

 The proof is obtained by induction. The base case $N=1$ follows from Corollary \ref{cor:BDP-sl}. 
 Now, suppose that the assertion holds for $N-1$, that is, for every $\epsilon>0$, there exists  {volume-preserving cocycle} $F^{(N-1)}$ over $f$ with $d_{\cC^0}(F^{(N-1)},F)<\epsilon$ and $(\boldsymbol{E},X,F^{(N-1)},f)$ contains $N-1$ periodic points $p_1,\ldots,p_{N-1}$ such that for every $i=1,2,\ldots,N-1$, $F^{(N-1)}_{n(p_i)} (p_i)=\id $. We aim to prove the statement for $N$. Fix $\epsilon>0$ and let $\delta_1:=\mathrm{min}\{\epsilon,\delta_0\}/3$.  
 {(Note that in this proof $F^{(i)}$ denotes the perturbations of cocycle $F$ and should not be mistaken with $F_i$ that corresponds to iterations of the cocycle.)}

 For every $i=1,2,\ldots,N-1$,  consider the map $\Theta_i :\sldr\to \sldr $ defined by 
 \[
 \Theta_i(H):=H . F^{(N-1)}(f^{n(p_i)-2}(p_i))\cdots F^{(N-1)}(f(p_i)) F^{(N-1)}(p_i).\] 
 
 It is easy to see that $\Theta_i$ is an open map. {Since we can approximate $\id$ matrix with hyperbolic matrices,} for every $1\leq i \leq N-1$, there exists $D_i\in \sldr$, $\delta_1$ close to  ${F}^{(N-1)}  (f^{n(p_i)-1}(p_i))$ such that $\Theta_i(D_i)$ is a hyperbolic matrix. Thus, applying Lemma \ref{lem:extension-bundle} with $K$ equal to the union of orbits of $p_i$'s ($1\leq i \leq N-1$), we find a {volume-preserving cocycle $\tilde{F}^{(N-1)}$ over $f$} with $d_{\cC^0}(\tilde{F}^{(N-1)},{F}^{(N-1)})<\delta_1$ such that for every $1\leq i\leq N-1$, 
\[
\tilde{F}^{(N-1)}(f^j(p_i))=\begin{cases}
    {F}^{(N-1)}(f^j(p_i)) & 0\leq j <n(p_i)-1\\
    D_i & j=n(p_i)-1
\end{cases}
\]

By this way of extension, we conclude that $\widetilde{F}^{(N-1)}_{n(p_i)}(p_i)=\Theta_i(D_i)$ is hyperbolic for every $i=1,\ldots,N-1$.   
 Then, for every $1\leq i\leq N-1$ consider the map 
 \(
 \Phi_i :\sldr^{n(p_i)-1}\to \sldr
 \) with 
 \[
 \Phi_i(A_{n(p_i)-1},\ldots,A_{1}):= A_{n(p_i)-1}\cdots A_2 A_1.
 \] 
 Clearly, $\Phi_i$ is an open map.  Now, since the set of hyperbolic matrices is open in $\sldr$ (see Remark \ref{rm:hyp-mat}), there exists  $\delta_2>0$ such that for every cocycle $H$ with $d(H,\widetilde{F}^{(N-1)}) < \delta_2$, all the linear maps $H_{n(p_i)}(p_i)$ for $i=1,\ldots,N-1$ are hyperbolic.
 
 Let $\delta_3<\mathrm{min}\{\delta_1,\delta_2\}$ be a positive number. Using Corollary \ref{cor:BDP-sl}, there is {volume-preserving cocycle $\hat{F}^{(N)}$ over $f$}  with $d_{\cC^0}(\hat{F}^{(N)},\widetilde{F}^{(N-1)})<\delta_3 $, and a periodic point $p_N\in X$ for $f$ such that $\hat{F}^{(N)}_{n(p_N)}(p_N)=\id $. Note that since $\delta_3<\delta_2$, ${\widetilde{F}}^{(N)}_{n(p_i)}(p_i)$ is hyperbolic for every $1\leq i \leq N-1$. This, in particular, implies that $p_N$ is not contained in the orbits of $p_i$'s for $i<N$. 
 \medskip

 Finally,  Lemma \ref{lem:extension-bundle}, for $K:= \mathcal{O}(p_1)\cup \mathcal{O}(p_2)\cup \cdots \cup  \mathcal{O}(p_N) $, implies that there exist a {volume-preserving cocycle ${F}^{(N)}$ over $f$} with $d_{\cC^0}({F}^{(N)},\hat{F}^{(N)})<\epsilon$, such that 
\begin{itemize}
    \item ${F}^{(N)}|_{\cO(p_N)}=\hat{F}^{(N)}|_{\cO(p_N)}$, 
    \item and for every $i=1,\ldots,N-1$, ${F}^{(N)}|_{\cO(p_i)}={F}^{(N-1)}|_{\cO(p_i)}$.
\end{itemize}
This completes the proof of Theorem \ref{thm:BDP-multi-orbits}.
 \end{proof}

\section{Proof of the dichotomy: domination vs. quasi-conformality}\label{sec:dichotomy}

Throughout this section, all the cocycles are $\sldr$ cocycles, defined on the trivial $d$-dimensional vector bundle over the base dynamics $\sigma:\Sigma_m\to \Sigma_m$. 

 This section is devoted to the proof of {Theorems \ref{thm:SL-B} and \ref{thm:lc-dichotomy}. These two theorems are direct corollaries of the following Theorem.}

\begin{theorem}\label{thm:dichotomy}
Assume that the cocycle $(\sigma,A)$ does not admit any dominated splitting. Then, for every $\epsilon>0$, there exists a {locally constant} cocycle $(\sigma,\Tilde{A})$ with $\|A-\Tilde{A}\|_{\mathcal{C}^0}<\epsilon$ such that  $(\sigma,\Tilde{A})$ stably admits quasi-conformal orbits in the $\cC^\alpha$ topology, that is for every $\alpha$-Hölder map $B$, sufficiently close to $\Tilde{A}$ in the $\cC^\alpha$ topology, $(\sigma,B)\in \qc$. 
 \end{theorem}

  This theorem provides the following dichotomy on approximating arbitrary  Hölder continuous cocycles over the shift.
  \begin{corollary}\label{cor:dichotomy}
  Let $(\sigma,A)$ be a Hölder continuous cocycle. Then, for any $\epsilon>0$, there exists a Hölder continuous cocycle $(\sigma,\Tilde{A})$ such that $\|A-\Tilde{A}\|_{\mathcal{C}^0}<\epsilon$ and 
  \begin{itemize}
      \item either  $(\sigma,\Tilde{A})$  admits a dominated splitting,
      \item or $(\sigma,\Tilde{A})$ is $\mathcal{C}^\alpha$-stably admits quasi-conformal orbits. 
  \end{itemize}
  \end{corollary}

  \subsection{Realization of arbitrary transitions}

In this subsection, we give some results on realizing a family of given matrices as the transitions of a perturbation of cocycles admitting periodic orbits with the identity action on the fiber. We first state the results for the cocycles acting with identity on the fiber at fixed points, then generalize them to those acting as the identity along periodic orbits (on the period).

\begin{proposition} \label{prop:realize-arbit-family-fix-pt}
Let $m\geq 2$ and $(\sigma,A)$ be a locally constant $\sldr$ cocycle over the shift map $\sigma:\Sigma_m\to \Sigma_m$ {such that  $A$ equals to the identity over the fixed points of $\sigma$}. Then, for every $\epsilon>0$, and every  $S_1,S_2,\ldots,S_{m-1}\in \sldr$, there exist
\begin{itemize}
    \item a positive integer $n$, 
    \item a locally constant $\sldr$ cocycle $(\Tilde{A},\sigma)$ depending on positions in $\ldbrack -n,n\rdbrack$ with $d_{\cC^0}(\Tilde{A},A)<\epsilon$,
    \item  $\uw_{i,j}\in T(\underline{p}_i, \underline{p}_j)$ for every $i\neq j$, where $\underline{p}_i:=(\underbrace{i,i,\ldots,i}_{2n+1})$.
\end{itemize}
  such that for every $1\leq i \leq m$ 
  \begin{equation}\label{eq:equality-muliset-}
  \{\tilde{A}_{[{2n+1;}\uw_{i,j}\rangle}:j\neq i\}=\{S_1,\ldots,S_{m-1}\}.
  \end{equation}
\end{proposition}
 \begin{remark}
     Note that equality in \eqref{eq:equality-muliset-} is equality of multi-sets. In other words, if some of $S_i$'s are equal, both sets are the same when elements are counted with multiplicities.
 \end{remark}

We need the following well-known fact for connected Lie groups, usually called fragmentation property.
  
\begin{lemma}{\cite[Corollary 2.9]{Kirillov_Book}}\label{lem:making-arbitrary-matrix}
    If $G$ is a connected Lie group and $U$ is a neighborhood of the identity element, then $U$ generates $G$. 
\end{lemma}

\begin{proof}[Proof of Proposition \ref{prop:realize-arbit-family-fix-pt}]
   Let $ l \in \mathbb{N} $ be such that the cocycle $ (\sigma,A) $ depends on the positions in $ \mathbb{I}_l $. We aim to provide $\tilde{A}$ as a locally constant cocycle that depends on positions in $ \mathbb{I}_n $ for a sufficiently large $ n > l $, which will be specified at the end of the proof. \\
   {
    Throughout the proof by  $A_{\transition{\uw_{i,j}}}$ we mean $A_{\transition{2n+1;\uw_{i,j}}}$ and remove $2n+1$ in the notation for simplicity.}
   
    For every $1\leq i \leq m$, let  $p_i:=(\ldots,i;i,i,\ldots)$ be the fixed point associated with $\underline{p}_i$.  For every distinct $i,j\in \{1,\ldots,m\}$  define 
    \begin{equation*}
   \uw_{i,j}:=   \big(\underbrace{i,i,\ldots,i}_{2n+1},\underbrace{j,j,\ldots,j}_{2n+1}\big). 
\end{equation*}
    Clearly, $\uw_{i,j}$ gives a transition from $\underline{p_i}$ to $\underline{p_j}$.
For every $1\leq i, j \leq m$ and $1\leq k \leq 2n$, let $\phi^{(k)}_{i,j}\in \sldr$ be such that 
\[\phi^{(k)}_{i,j}:=A|_{\cylinder[-n;(\underbrace{i,\ldots,i}_{2n+1-k},\underbrace{j,\ldots,j}_{k})]}.\]
\begin{figure}[t]
    \centering    \includegraphics[width=.7\textwidth]{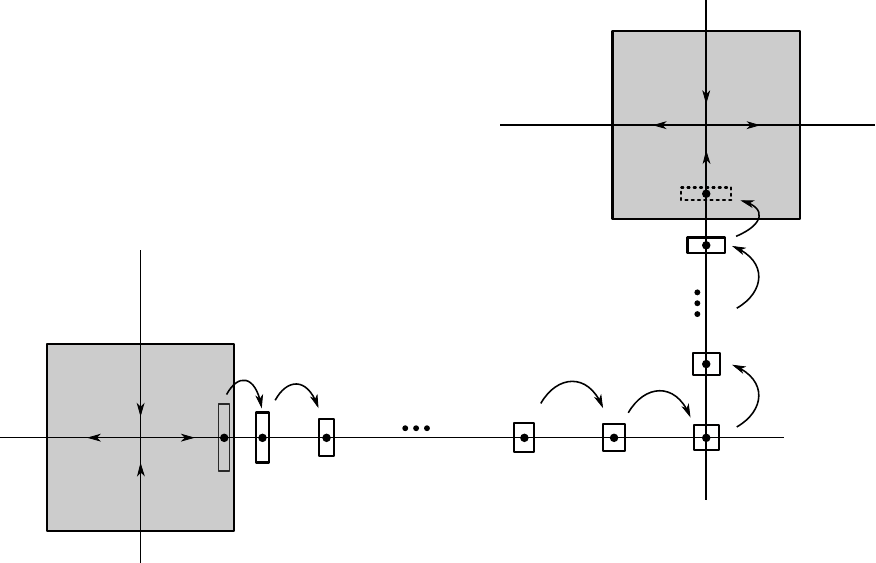}
    \caption{A transition following a homoclinic point. On grey regions, the cocycle is identity. The dashed rectangle is zoomed in the next figure.}
    \label{fig:transition}
\end{figure}

Since $A(p_i)=A(p_j)=\id $, and $A$ depends on positions in $\cA^{\ldbrack -l,l\rdbrack}$, for every $k$ with $k\leq n-l$ or $k\geq n+1+l$, $\phi^{(k)}_{i,j}=\id $. Hence, 
\begin{equation*}A_{\transition{\uw_{i,j}}}=\underbrace{\phi^{(2n)}_{i,j}\cdots \phi^{(n+l+1)}_{i,j}}_{\id } \phi^{(n+l)}_{i,j}\cdots \phi^{(n-l+1)}_{i,j} \underbrace{\phi^{(n-l)}_{i,j} \cdots \phi^{(1)}_{i,j}}_{\id }=\phi^{(n+l)}_{i,j}\cdots \phi^{(n-l+1)}_{i,j},
\end{equation*}

For every $1\leq i \leq m$, let $\{D_{i,j}:1\leq j \leq m, j\neq i\}$ be an arbitrary permutation of $S_1,\ldots,S_{m-1}$. 
One can apply Lemma \ref{lem:making-arbitrary-matrix} to each $(A_{\transition{\uw_{i,j}}})^{-1}D_{i,j}$ ($1\leq i ,j\leq m$, $i\neq j$), and obtain $n_{i,j}\in \N$ and an $n_{i,j}$-tuple $(D^{(1)}_{i,j},\ldots,D^{(n_{i,j})}_{i,j})$ of elements of $\sldr$, $\epsilon$-close to $\id$ such that for every $1\leq i,j\leq m$ with $i\neq j$, and $1\leq s \leq n_{i,j}$
\begin{equation}\label{eq:prop-D_ij}
(A_{\transition{\uw_{i,j}}})^{-1}D_{i,j}=D^{(n_{i,j})}_{i,j}D^{(n_{i,j}-1)}_{i,j}\cdots D^{(1)}_{i,j}, \quad \|D_{i,j}^{(s)}-\id\|<\epsilon.
\end{equation}
Note that without loss of generality, one may assume that all the numbers $n_{i,j}$ have equal value, say $n_0$. Indeed, for given $n_0>\max\{n_{i,j}\}$,  one can define $D^{(s)}_{i,j}=\id$ for $n_{i,j}<s\leq n_0$. By this definition, all the properties listed in \eqref{eq:prop-D_ij} hold.
\medskip
Now, we take $n:=n_0+l$ and define $\tilde{A}:\Sigma_m\to \sldr$ depending on positions in $\cA^{\ldbrack -n,n\rdbrack}$ such that for every $1\leq i,j\leq m$ and $1\leq k \leq n-l$, 
\[\tilde{A}|_{\cylinder[-n;\underbrace{i,\ldots,i}_{2n+1-k},\underbrace{j,\ldots,j}_{k}]}:=D^{(k)}_{i,j},\]
and for all other cylinders $\cylinder[-n;\ua]$ with $\ua\in \cA^{\ldbrack -n,n\rdbrack}$, $\tilde{A}|_{\cylinder[-n;\ua]}:=A|_{\cylinder[-n;\ua]}$.
\medskip
We claim that $\tilde{A}$ satisfies all the required conditions. First note that since all the matrices $D_{i,j}^{(s)}$ are $\epsilon$-close to the identity, $d_{\cC^0}(A,\tilde{A})<\epsilon$. On the other hand, for every $1\leq i,j\leq m$ with $i\neq j$, 
\begin{align*}\tilde{A}_{\transition{\uw_{i,j}}}=\underbrace{\phi^{(2n)}_{i,j}\cdots \phi^{(n+l+1)}_{i,j}}_{\id } \underbrace{\phi^{(n+l)}_{i,j}\cdots \phi^{(n-l+1)}_{i,j} }_{{A}_{\transition{\uw_{i,j}}}}\underbrace{D^{(n-l)}_{i,j} \cdots D^{(1)}_{i,j}}_{({A}_{\transition{\uw_{i,j}}})^{-1}D_{i,j}}=D_{i,j}.
\end{align*} 
Hence, by definition of $D_{i,j}$'s, for every $1\leq i\leq m$, 
\[\{\tilde{A}_{\transition{\uw_{i,j}}} :  j\neq i \}=\{D_{i,j} :  j\neq i \}=\{S_1,\ldots,S_{m-1}\}.\]
 This completes the proof of the proposition.
\begin{figure}[t]
    \centering 
    \includegraphics[width=.55\textwidth]{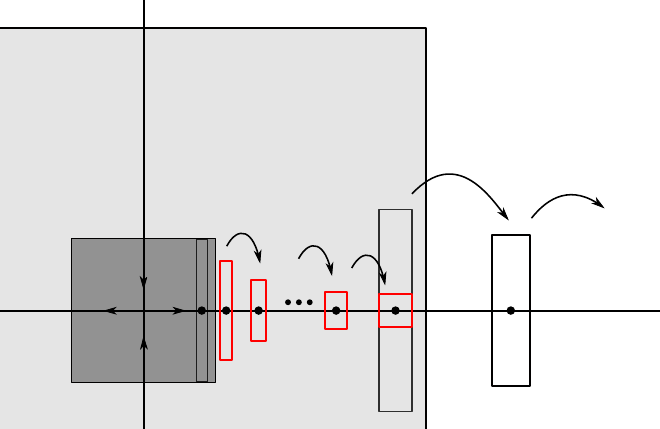}
    \caption{Realization of arbitrary transition matrix in an identity region. The perturbed cocycle in each red region is a matrix near the identity, whose product can be an arbitrary element of $\sldr$, while on the dark grey region, it remains identity.}
    \label{fig:transition-arbitrary-realization}
\end{figure}
\end{proof}

\begin{proposition} \label{prop:realize-arbit-family-per-pt}
Let $(\sigma, A)$ be a locally constant $\sldr$ cocycle over
{   the full shift on finitely many symbols
    such that there exist periodic points $p_1,p_2,\ldots,p_m$ {with disjoint orbits such that for every $i=1,\ldots,m$}, $A_{n(p_i)}(p_i)=\id$.}

Then, for every $\epsilon>0$, and every  $S_1,S_2,\ldots,S_{m-1}\in \sldr$, there exist
\begin{itemize}
    \item positive integers $n,k$, 
    \item a locally constant $\sldr$ cocycle $(\sigma,\Tilde{A})$ depending on positions in {$\ldbrack -n, n\rdbrack$} with $d_{\cC^0}(\Tilde{A},A)<\epsilon$,
    \item $\ua_1,\ldots,\ua_m\in \cA^{\ldbrack -n, n\rdbrack}$ so that $p_i$ is the periodic point associated to $\ua_i$.
    \item  $\uw_{i,j}\in T(\ua_i, \ua_j)$ for every $i\neq j$
\end{itemize}
  such that for every $1\leq i \leq m$ 
  \begin{equation}\label{eq:equality-muliset}
      \{\tilde{A}_{\transition{k;\uw_{i,j}}} :  j\neq i\}=\{S_1,\ldots,S_{m-1}\}.
  \end{equation}
\end{proposition}
\begin{proof}
    Take $N:=|\ua_1|.|\ua_2|\cdots |\ua_m|$. 
    For every $i=1,\ldots,m$, let $$\ub_i:=\underbrace{\ua_i\,\ua_i\,\cdots \,\ua_i}_{\frac{1}{|\ua_i|}N}\in \cA^{N}.$$ Then, $\underline{p}_i$ is the periodic point associated to $\ub_i$. Note that for every $i=1,\ldots,m$, $A_N(p_i)=\id$, as $N$ is divisible by $n(p_i)$.  Define the locally constant cocycle  $B:\Sigma_m\to \sldr$ by 
    \[B(\ldots,j_{-1};j_0,j_1,\ldots):=A_N(\cdots\ub_{j_{-1}};\ub_{j_{0}}\,\ub_{j_{1}}\cdots ).\]
Note that since $A$ is locally constant, there exists $T \in \N$ such that $A$ depends on positions in $\ldbrack -T N,(T+1)N-1\rdbrack$. Therefore, the cocycle $B$ depends on positions in $\ldbrack -T, T\rdbrack$. Moreover, denoting $p'_i=(\ldots i; i, i, i,\ldots)$, we have $B(p'_i)=A_N(p_i)=\id$ for every $i=1,\ldots,m$.

Now, $B$ satisfies all the conditions of Proposition \ref{prop:realize-arbit-family-fix-pt} and by this proposition, there exists $n'\in \N$, a continuous cocycle $\tilde{B}:\Sigma_m\to \sldr$ with $d_{\cC^0}(\tilde{B},B)<\epsilon\|A\|^{-Nd}$  depending on positions in $\ldbrack -n',n'\rdbrack$, and transitions $\uw'_{i,j}$ from $(i,i,\ldots,i)\in \cA^{2n'+1}$ to $(j,j,\ldots,j)\in \cA^{2n'+1}$ such that
    \begin{equation}\label{eq:desired-transition-B}
        \{\tilde{B}_{[2n'+1;\uw'_{i,j}\rangle}:j\neq i\}=\{S_1,\ldots,S_{m-1}\}, \quad \text{for every }i=1,\ldots,m.
    \end{equation}
    Now, we define the cocycle $\tilde{A}$. 
    \begin{itemize}
        \item For any point $x$ in a cylinder of the form $\cylinder[-Nn';\ua_{j_{-n'}}\ldots\ua_{j_0}\ldots\ua_{j_{n'}}]$, let $\tilde{x}$  be a point in $\cylinder[-n';(j_{-n'},\ldots; j_{0},\ldots,j_{n'})]$ and set 
        \[\tilde{A}(x):=\tilde{B}(\tilde{x})B(\tilde{x})^{-1}A(x).\]
        Note that the definition is independent of the choice of $$\tilde{x}\in \cylinder[-n';(j_{-n'},\ldots; j_{0},\ldots,j_{n'})]$$ as both $B,\tilde{B}$ are constant on this cylinder. 
        \item For all other point $x$, define $\tilde{A}(x):=A(x)$. 
    \end{itemize}
    We claim that this definition satisfies all the required properties for $\tilde{A}$. 
    
    \medskip
    \noindent
        $\bullet$ {\bf $\tilde{A}$ is locally constant.} It is clear from the definition that $\tilde{A}$ is a locally constant cocycle depending on positions in $\ldbrack -Nn',N(n'+1)\rdbrack$.

    \medskip
    \noindent
        $\bullet$ {\bf $d_{\cC^0}(A,\tilde{A})<\epsilon$.} We only need to check that for any element $x$ of a cylinder of the form $\cylinder[-Nn';\ua_{j_{-n'}}\ldots\ua_{j_0}\ldots\ua_{j_{n'}}]$, $\|A(x)-\tilde{A}(x)\|<\epsilon$ (otherwise $A=\tilde{A}$). For such $x$, \[\tilde{A}(x)-A(x)=\tilde{B}(\tilde{x})B(\tilde{x})^{-1}A(x)-A(x)=\big(\tilde{B}(\tilde{x})-B(\tilde{x})\big)B(\tilde{x})^{-1}A(x).\]
        Hence, as we know $d_{\cC^0}(B,\tilde{B})<\epsilon \|A\|^{-Nd}$, 
        \begin{align*}
            \|\tilde{A}(x)-A(x)\| \leq \underbrace{\|\tilde{B}(\tilde{x})-B(\tilde{x})\|}_{\leq \epsilon \|A\|^{-Nd}}.\underbrace{\|B(\tilde{x})^{-1}\|}_{\leq \|A\|^{N(d-1)}}.\|A(x)\|<\epsilon,
        \end{align*}
        where the upper bound for $\|B(\tilde{x})^{-1}\|$ comes from the fact that for every $D\in \sldr$, $\|D^{-1}\|\leq \|D\|^{d-1}$. 
        
    \medskip
    \noindent
        $\bullet$ 
        {\bf Transitions.} If $\uw'_{i,j}=(j_1,\ldots,j_n)$, then define 
         \begin{equation*}
        \uw_{i,j}:=\ua_{j_1}\,\ua_{j_2}\,\ldots\,\ua_{j_n}.
        \end{equation*}
        Then, $\tilde{A}|_{[k;\uw_{i,j}\rangle}=\tilde{B}|_{[2n'+1;\uw'_{i,j}\rangle}$, where $k=(2n'+1)N$.  This completes proof due to \eqref{eq:desired-transition-B}.
    \end{proof}

\begin{corollary}\label{cor:make-gen-cov}
   Let $\ell>d^2$ and $(\sigma,A)$ be a locally constant $\sldr$ cocycle over
   the full shift on finitely many symbols
    such that there exist periodic points $p_1,p_2,\ldots,p_\ell$ {with disjoint orbits such that} $A_{n(p_i)}(p_i)=\id$, for every $i=1,\ldots,\ell$. Then, for every $\epsilon>0$, there exists a locally constant cocycle $(\sigma,\Tilde{A})$ with  $d_{\cC^0}(\Tilde{A},A)<\epsilon$ satisfying  covering condition.
\end{corollary}

\begin{proof}
It suffices to take $S_1,S_2,\ldots,S_{\ell-1}$ in Proposition \ref{prop:realize-arbit-family-per-pt} as the family given by Lemma \ref{lem:sldr-covering}. We know that for this family, there exists a non-empty open subset $\cU\subset \sldr$ with compact closure such that 
\[\overline{\cU}\subset \bigcup_{i=1}^\ell S_i^{-1}\cU.\]
Then, according to Proposition \ref{prop:realize-arbit-family-per-pt}, one can perturb $A$ in the $\cC^0$ topology to get a cocycle $\tilde{A}$ such that $(\sigma,\tilde{A})$ satisfies the covering condition.  
\end{proof}

\subsection{Proof of Theorem \ref{thm:dichotomy}}

This subsection is dedicated to the proof of Theorem \ref{thm:dichotomy}.

We begin with the following lemma, which asserts that arbitrary $\cC^0$-close to every continuous cocycle, it is possible to find a locally constant cocycle whose values equal the original cocycle at a predetermined finite set of points. {The proof is simple and is left to the reader. 
}

\begin{lemma}\label{lemma-perturb-to-loc}
Let $(\sigma,A)$ be a continuous cocycle over the shift map $\sigma:\Sigma_m\to \Sigma_m$ and $K\subset \Sigma_m$ be a finite set. Then, for every $\epsilon>0$, there exists a locally constant map $B:\Sigma_m\to \sldr$ such that $d_{\cC^0}(A,B)<\epsilon$ and $B|_K=A|_K$. 
\end{lemma}

Now, we have all the ingredients to provide the proof of Theorem \ref{thm:dichotomy}. 

\begin{proof}[Proof of Theorem \ref{thm:dichotomy}]
The proof is based on the combination of Corollary \ref{cor:stably no dom}, Theorem \ref{thm:BDP-multi-orbits}, Lemma \ref{lemma-perturb-to-loc},  Corollary \ref{cor:make-gen-cov} and Theorem \ref{thm:SFT-covering}.

{According to Corollary \ref{cor:stably no dom}, we can assume that $(\sigma,A)$ stably does not admit dominated splitting.} Now, it follows from Theorem \ref{thm:BDP-multi-orbits} with $N=d^2$ that there exists a cocycle $A^{(1)}$ with $d_{\cC^0}(A^{(1)},A)<\epsilon/3$ and periodic points $p_1,\ldots,p_N$ (with disjoint orbits) such that  for $i=1,\ldots,N$ that $A^{(1)}_{n(p_i)}(p_i)=\id$.

Then, applying Lemma \ref{lemma-perturb-to-loc} to $A^{(1)}$ for $K$ equal to union of orbits of $p_i$'s, we get a locally constant cocycle $A^{(2)}$ with $d_{\cC^0}(A^{(2)},A^{(1)})<\epsilon/3$ such that for $i=1,\ldots,N$ that $A^{(2)}_{n(p_i)}(p_i)=\id$.

Finally, by Corollary \ref{cor:make-gen-cov},  there exists a {locally constant} cocycle $A^{(3)}$ with $$d_{\cC^0}(A^{(3)},A^{(2)})<\epsilon/3$$ satisfying the covering condition with respect to some open set $\cU\subset \sldr$ with compact closure. Hence, $d_{\cC^0}(A^{(3)},A)<\epsilon$ and  by Theorem \ref{thm:SFT-covering}, for every $\alpha>0$, every cocycle sufficiently $\cC^\alpha$-close to $(\sigma,A^{(3)})$ admits quasi-conformal orbits. 
\end{proof}

 \section{Further comments and questions}\label{sec:questions}
 
\subsection{Iterated Function Systems}\label{subsec:question-IFS}

It is natural to expect a dichotomy similar to Theorem \ref{thm:SL-B} holds for Iterated Function Systems (IFS) generated by generic pairs (or \(n\)-tuples) of matrices in $\sldr$. This can be seen as the first step toward understanding the behavior of generic cocycles with high regularity.

\begin{question}\label{ques:generic-IFS}
Let \(n \geq 2\), and \((H_1, \dots, H_n)\) be a generic \(n\)-tuple of matrices in \(\sldr\).
    Assume that the IFS generated by \(H_1, \dots, H_n\) has no bounded orbit branch. Does the IFS exhibit uniform exponential growth? Does it admit a dominated splitting?
\end{question}
Recall that the IFS generated by \(H_1, \dots, H_n\) does not have a bounded orbit branch if for every $j_1,j_2,\ldots\in \{1,2,\ldots,n\}$, the sequence $\|H_{j_k}\cdots H_{j_1}\|$  diverges to infinity. Moreover,  uniform exponential growth means that the matrices in the semigroup generated by \(H_1, \dots, H_n\)  grow uniformly exponentially with word length. More precisely, there exist $C>0,\lambda>1$ such that $j_1,j_2,\ldots\in \{1,2,\ldots,n\}$, the sequence $\|H_{j_k}\cdots H_{j_1}\|>C\lambda^k$. 

 We remark that Question \ref{ques:generic-IFS} for $d=2$ has a positive answer due to the work of \cite{Avila_Bochi_Yoccoz} (see also \cite{Christodoulou}). They show that if the semigroup generated by $n$-tuple $(H_1,H_2,\ldots,H_n)\in \mathrm{SL}(2,\R)^n$ does not contain any elliptic elements, after a sufficiently small perturbation of the generators, the action of matrices on $\R^2$, leaves a multi-cone invariant. By a multi-cone, we mean a union of finitely many closed cones in $\R^2$ whose mutual intersection is only the origin. The existence of an invariant multi-cone implies the uniform exponential growth of norms for the cocycle. {For the case \(d=3\), see \cite{Potrie-et-al-SL3R}, where related questions about for finitely generated subgroups of \(\mathrm{SL}(3,\R)\) are investigated.}

 \subsection{Stability in low regularity
 }\label{subsec:question-C0-stable}
The analysis in the proof of our main theorems regarding \(\cC^\alpha\) stability of bounded behavior (cf. Section \ref{sec:stability-qc}) relied crucially on Hölder estimates. Thus, asking questions regarding the prevalence of cocycles with stable quasi-conformality in  \(\cC^0\) regularity is natural.

\begin{question}
    Let $d>2$. Does there exist an $\sldr$ cocycle which \(C^0\) stably exhibits bounded orbits?
\end{question}

{Even though the answer to this question may depend on the base dynamics, we expect a negative answer for any base if $d>2$. Corollary \ref{cor:minimal} provides a negative answer to this question for minimal bases. {After the first version of this paper appeared on arXiv, Bochi \cite{Bochi-no-stable-CQ} announced a negative answer to this question for cocycles over homeomorphisms of zero-dimensional topological spaces having finitely many periodic points of any period. Indeed, he showed that for such base systems, generic $\cC^0$ cocycles admit no quasi-conformal orbits. This, in particular, shows that the regularity in Theorem \ref{thm:GL-bd-example} is optimal.}  
Note that for $d=2$, the question has a positive answer when the base dynamics contain periodic points, as elliptic elements over periodic points (on the period) cannot be destroyed with small $\mathcal{C}^0$ perturbations.
}

\subsection{Finer classification of general linear  cocycles} \label{subsec:question-GL}

We expect further refinement of Theorem \ref{thm:GL-QC} should be true.

Let $\cV^+$ be the set of all uniformly volume expanding $\GL(d,\R)$ cocycles over SFT, that is, the set of all cocycles  $(\sigma,A)$ such that $\inf_x |\det(A_n(x))|>1$ for some $n\in\N$.

Given $(\sigma,A) \in \cV^+$, by {Theorem \ref{thm:GL-QC},  either it admits} a dominated splitting or $\cC^0$-approximated by an $\cC^\alpha$ open set of cocycles admitting $\kappa$-quasi-conformal orbits for some $\kappa>1$. Moreover, such orbits are uniformly expanding.

The same dichotomy holds in the set of all uniformly volume contracting $\GL(d,\R)$ cocycles over SFT, denoted by $\cV^-$, that is, the set of all cocycles  $(\sigma,A)$ such that $\sup_x |\det(A_n(x))|<1$ for some $n\in \N$.

We expect a positive answer to the following question. 
 \begin{question}
\label{question:GL}
Let $\cB$ be the set of $\GL(d, \R)$ cocycles over subshift of finite type that admit bounded orbits. How big is the $\cC^\alpha$-interior of $\cB$? Is it $\cC^0$-dense in the complement of $\cD\cup\cV^-\cup\cV^+$?
\end{question}
Whenever a cocycle is far from being either uniformly volume contracting or uniformly volume expanding, there are simultaneously volume contracting and volume expanding orbits. In this case, if, in addition, the cocycle does not admit any dominated splitting, then by \cite{BDP} arguments (i.e. Proposition \ref{prop:BDP-dense-dichotomy} and proof of Theorem \ref{thm:BDP-multi-orbits}) we can show that the cocycle is approximated by another cocycle with the same properties with several periodic points which are homothety at their periods. 
This is not enough to address Question \ref{question:GL}, and one must show that the cocycle is approximated by another cocycle with the same properties and several periodic points, which are \emph{identity} at their periods.

\appendix

\section{}

\subsection{Generic properties of linear cocycles} 
\label{sec:generic SLd}
{In this section, we explain how a generic cocycle over the shift map behaves at periodic orbits. Indeed, we will prove that the iterations of a generic cocycle are unbounded on every periodic fiber. }

    \begin{lemma}\label{lem:generic-pair-hyperbolic}
For every $m,d\in \N$ with $d>2$, for a generic $m$-tuple  $(A_1,A_2,\ldots,A_m)$ in $ \sldr^m$, all the elements of the semigroup generated by $A_1,\ldots,A_m$ are hyperbolic. 
\end{lemma}

 \begin{proof}
     {Fix a periodic word $\uw = (w_1, \dots, w_k) \in \cA^{\mathrm{fin}}$, and consider the product map  
\[
P_{\uw}[W] = W_{w_k} \cdots W_{w_1}, \quad [W] = (W_1, \dots, W_m) \in \SL(d,\mathbb{R})^m.
\]  
We claim the set  
\(
\cH_{\uw} := \{ [W] \in \SL(d,\mathbb{R})^m : P_{\uw}[W] \text{ is hyperbolic} \}
\)
is open and dense in $\SL(d,\mathbb{R})^m$.  
Note that $P_{\uw}[W]$ is hyperbolic if none of its eigenvalues lie on the unit circle. Now, define  
\[
h_{\uw}[W] := \prod_{i,j=1}^d \left( \lambda^{(i)}_{\uw}[W] \lambda^{(j)}_{\uw}[W] - 1 \right).
\]  
It is a symmetric polynomial of the eigenvalues of $P_{\uw}[W]$. Hence, by the fundamental theorem of symmetric polynomials, $h_{\uw}[W]$ is a polynomial in the entries of $W_1, \dots, W_m$. One can easily show that its zero set corresponds to non-hyperbolic products.  \\
This zero set is a proper real-algebraic subset of $\SL(d,\mathbb{R})^m$, because one can explicitly choose diagonal hyperbolic matrices (e.g., $W_i = \mathrm{diag}(2, 1/2, 1, \dots, 1)$) such that $P_{\uw}[W]$ is hyperbolic and $h_{\uw}[W] \neq 0$. Therefore, $\cH_{\uw}$ is open and dense.  
Finally, since $\cA^{\mathrm{fin}}$ is countable, and for every $\uw \in \cA^{\mathrm{fin}}$, $\cH_{\uw}$ is open and dense, the set 
$\bigcap\limits_{\uw\in \cA^{\mathrm{fin}}} \cH_{\uw}$ is a generic set. }
 \end{proof}

\begin{remark} \label{rm:hyp-mat}{
For $d > 2$, hyperbolic matrices form an open and dense subset of $\SL(d,\mathbb{R})$. When $d = 2$, they are open but not dense; both hyperbolic and elliptic matrices are open, and their union is dense. Also, for $d > 2$, the proof Lemma \ref{lem:generic-pair-hyperbolic} implies  that, for almost every $(A_1, \ldots, A_m)$ in $ \SL(d,\mathbb{R})^m$ (with the respect to the product of Haar measures), the semigroup $\langle A_1, \ldots, A_m\rangle^+$ consists only of hyperbolic matrices.
}\end{remark}

 {The next proposition shows that in all regularities a typical $\sldr$ cocycle has only hyperbolic periodic fibers if $d>2$.}

\begin{proposition}\label{prop:generic hyp periodic}
Let $X$ be a compact metric space and $f:X\to X$ be a homeomorphism with countably many periodic points. Also, let $d> 2$ and $r\in [0,1]$. Then,  every {
$\cC^r$}
 generic  $\sldr$ cocycle  $(f,A)$ satisfies the following:\\
 For every periodic point $p$ for $f$ of period $n(p)$, the matrix $A_{n(p)}(p)$ is hyperbolic. 
\end{proposition}

\begin{proof}
    {
For periodic point $p$, the set of cocycles $(f,A)$ for which $A_{n(p)}(p)$ is hyperbolic is open in the $\cC^0$ topology (hence open in the $\cC^r$ topology for every $r\in [0,1]$), since hyperbolic matrices are open in $\SL(d,\mathbb{R})$.\\
To show density in the $\cC^r$ topology, assume $A_{n(p)}(p)$ is not hyperbolic. Since $d>2$, one can easily construct a smooth path $H_t \in \SL(d,\mathbb{R})$ with $H_0 = \id$ and $A_{n(p)}(p) H_t$ hyperbolic for all $t > 0$ (see \citep[Lemma 2.20]{Reshadat-thesis} for the details).
 Let $\psi : X \to [0,1]$ be a $\cC^r$ function with $\psi(p) = 1$ and $\psi \equiv 0$ outside a  $\epsilon$-neighborhood of $p$. Define
\(
A^{(t)}(x) := A(x) H_{t\psi(x)}.
\)
Then  for small $t > 0$, $A^{(t)}$ is $\cC^r$-close to $A$ and agrees with $A$ outside $\epsilon$-neighborhood of $p$. Also, $\big(A^{(t)}\big)_{n(p)}(p) = A_{n(p)}(p) H_t$, which is hyperbolic.
Finally, since the periodic points of $f$ are countable, the conclusion follows by the Baire category theorem.}
\end{proof}

\subsection{Extending a cocycle}\label{sec:extension}

This subsection concerns results on extending a linear cocycle defined on a proper subset to the ambient space $X$. 
   \begin{lemma}\label{lem:extension}
Let $X$ be a normal topological space. Then, for every $\epsilon>0$, there exists $\delta>0$, such that for every compact subset $K\subset X$, whenever  $A:X\to \sldr$, $\tilde{A}_K:K\to \sldr$ are bounded continuous maps with  $\sup_{x\in K}\|A(x)-\tilde{A}_K(x)\|<\delta$, there exists continuous map $\tilde{A}:X\to \sldr$ satisfying \[\left.\Tilde{A}\right|_K=\Tilde{A}_K~~ \text{and} ~~ \sup_{x\in X}\|A(x)-\Tilde{A}(x)\|<\epsilon.\]
 \end{lemma}

\begin{proof}{
For $1 \leq i, j \leq d$, write $A_{i,j}$ and $\tilde{A}_{K,i,j}$ for the entries of $A$ and $\tilde{A}_K$, respectively. Define $H_{K,i,j} := \tilde{A}_{K,i,j} - A_{i,j}|_K$, which is continuous on the compact set $K$. By Tietze’s extension theorem, each $H_{K,i,j}$ extends to a continuous function $H_{i,j}:X \to \R$ with $\sup_{x\in X}  |H_{i,j}(x)|=\sup_{x\in K} |H_{K,i,j} (x)|$.\\
Let $B_{i,j} := A_{i,j} + H_{i,j}$, and define the matrix-valued map $B = (B_{i,j}): X \to \mathrm{M}_{d \times d}(\R)$. Then $B$ is continuous, $B|_K = \tilde{A}_K$, and for sufficiently small $\delta>0$, if $\sup_{x\in K}\|A(x) - \tilde{A}_K(x)\|< \delta$, then for every $x\in X$, $\|A(x) - B(x)\| < \epsilon/2$ and $|\det B(x)-1|<\epsilon/2$. Note that $\det B(x)$ is not necessarily equal to 1, so we can normalize $B$ by defining 
\[
\tilde{A}(x) := \left( \det B(x) \right)^{-1/d} B(x)\in \sldr,
\]
Now, $\tilde{A}: X \to \sldr$ is continuous and if $\delta$ is sufficiently small (depending on $\epsilon$ and $\|A\|$), then $\sup_{x\in X}\|\tilde{A}(x) - A(x)\| < \epsilon$. 
On the other hand, $\tilde{A}$ coincides with $A$ on $K$, as for every $x\in K$, $B(x)=A(x)$ and $\det B(x)=1$. 
}\end{proof}

\begin{lemma}\label{lem:extension-bundle}
      Let $(\boldsymbol{E},X,F,f)$ be a continuous volume-preserving cocycle defined on a normal topological space $X$. Then, for every $\epsilon>0$, there exists $\delta>0$ such that for every   finite set $K\subset X$, invariant under $f$ and for every {volume-preserving} $\delta$-perturbation (in the $\cC^0$ topology) of $F$ on $K$, denoted by $\tilde{F}_K$, there is an $\epsilon$-perturbation (in the $\cC^0$ topology) $\tilde{F}$ of $F$ such that $\tilde{F}_K=\tilde{F}\vert_K$.
 \end{lemma}
\begin{proof}{
Let $K = \{x_1, \dots, x_k\}$ and choose 
small numbers $r>r'>0$ such that if $U_i,U'_i$ denotes the balls of radius of $r,r'$ around $x_i$, respectively, then 
$f(U'_i) \subset \bigcup_j U_j$ and the bundle $\boldsymbol{E}$ is trivialized over $U_i$ via homeomorphisms $\phi_i$. Define $U := \bigcup_i U_i$, $U' := \bigcup_i U'_i$. Define $\phi: U' \times \mathbb{R}^d \to \pi^{-1}(U')$ so that $\phi=\phi_i$ on $U'_i\times \R^d$. 
Using this trivialization, $F$ corresponds to a $\sldr$ cocycle $A: \overline{U'} \to \mathrm{SL}(d,\mathbb{R})$ via
\begin{equation}\label{eq:def-A}
    \phi^{-1} \circ F \circ \phi(x,v) = (f(x), A(x)v).
\end{equation}
Given $\epsilon > 0$, by Lemma~\ref{lem:extension}, there exists $\delta > 0$ such that  any $\delta$-perturbation $\tilde{A}_K$ of $A$ on $K \cup \partial U'$ extends to a continuous map $\tilde{A}: \overline{U'} \to \mathrm{SL}(d,\mathbb{R})$ with $\sup_{x\in \overline{U'}}\|\tilde{A}(x) - A(x)\| < \epsilon$.\\
Now, given a $\delta$-perturbation $\tilde{F}_K$ of $F$ on $K$, define $\tilde{A}_K$ on $K$ (via $\phi$) similar to \eqref{eq:def-A} and set $\tilde{A}_K = A$ on $\partial U'$. Then, apply
 Lemma \ref{lem:extension} to an extension $\Tilde{A}:\overline{U'} \mapsto \sldr$ of $\tilde{A}_K$ and define $\Tilde{F}$ by
\[
\tilde{F}(x,v) :=
\begin{cases}
(f(x), \phi \tilde{A} \phi^{-1}(x,v)), & \text{if } x \in \pi^{-1}(\overline{U'}), \\
F(x,v), & \text{otherwise}.
\end{cases}
\]
It is straightforward to verify that $\tilde{F}$ satisfies all the required properties, completing the proof.
}\end{proof}

\providecommand{\href}[2]{#2}

\end{document}